\title{ Extremal functions for the sharp $L^2-$ Nash inequality}
\author {Emmanuel Humbert\ \\  Institut \'Elie Cartan, Universit\'e de Nancy
  1, BP 239 \\ 54506 Vandoeuvre-L\`es-Nancy Cedex, FRANCE \\
Email : humbert@iecn.u-nancy.fr  }

\date{  }
\documentclass[10pt]{article}
\newtheorem{theorem}{Theorem}

\newtheorem{step}{Step}
\newtheorem{cor}{Corollary}
\font\tenbb=msbm10

\def\rR{\hbox{\tenbb R}}
\def\nN{\hbox{\tenbb N}}

\begin{document}

%%%%%%%%%%%%%%%%%%%%%%%%%%%%%%%%%%%%%%%%%%%%%%%%%%%%%%%%%%%%%%%%
\maketitle
%%%%%%%%%%%%%%%%%%%%%%%%%%%%%%%%%%%%%%%%%%%%%%%%%%%%%%%%%%%%%%%%
\begin{abstract}
We give geometrical conditions under which there exist
extremal functions for the sharp $L^2$-Nash inequality.
\end{abstract}

\section{Introduction}
This paper is in the spirit of several works on best constants problems
in Sobolev type inequalities. A general reference on this subject is
the recent book of Hebey \cite{heb}. These questions have many interests.
At first, they are at the origin of the resolution of famous geometrical 
problems as Yamabe problem. More generally, 
they show how geometry and analysis 
interact on Riemannian manifolds and lead to the developpement of 
interesting analytic 
methods. This article is devoted to the  existence of 
extremal functions for the optimal $L^2$-Nash inequality and 
follows another paper \cite{h} in which we proved the existence of a
second best constant in the $L^2$-Nash inequality. Obviously, 
finding extremal functions 
is interesting from PDEs'  point of view. The proof 
we give here may appear
very technical. Nevertheless, its interest lies in 
the analytic methods it gives, for example on what concerns the study
of concentration phenomenons. Moreover, extremal functions have their 
own interests because they give informations on best constants. For example, 
the existence of extremal functions for the circle $S^1$ gives an explicit
inequality on $S^1$ (see \cite{h}).
 
In this paper, 
we let $(M,g)$ be a smooth compact Riemannian $n$-manifold. We consider
the following inequality : for  $u \in C^{\infty}(M)$,
$$(\int_{M}u^2dv_g)^{1+\frac{2}{n}}\leq (A\int_{M}{ \mid \nabla u \mid}_g^2 dv_{g} +
B\int_{M}u^2dv_g)(\int_{M}\mid u \mid dv_g)^{\frac{4}{n}} \eqno { N (
  A ,B ) (u)}$$
We say that $N(A,B)$ $is$ $valid$ if $N(A,B)(u)$ is true for
all $u \in C^{\infty}(M)$. In the following, we refer to this
inequality as the $L^2$-Nash inequality. Let now 
$$A_0=\inf \{ A>0 | \hbox{ there exists } B>0 \hbox{ s.t. }
 N(A,B) \hbox{ is
  valid } \}$$
It was shown in \cite{cl} that 
$$A_0=A_0(n)=\frac{{(n+2)}^{\frac{n+2}{n}}}{2^{\frac{2}{n}}n \lambda_1(\mathcal{B})
{\mid \mathcal{B} \mid}^{\frac{2}{n}}}$$
where
$\mid \mathcal{B} \mid$  is the volume of the unit ball $\mathcal{B}$ in $\rR^n$,
${\lambda}_1$ is the first nonzero Neumann 
eigenvalue of the Laplacian for radial
functions on
$\mathcal{B}$ and $Vol(M)$ is the volume of $(M,g)$.
Then, it was shown in \cite{h} 
that there exists $B>0$ such that the sharp $N(A_0(n),B)$ is
valid. Another form of sharp inequality is in Druet-Hebey-Vaugon
\cite{dhv}.
Let now
$$B_0= \inf \{    B \in \rR \hbox{ s.t. } N(A_0(n),B) \hbox{  is valid }\}$$
It was also proved in \cite{h} that for any smooth compact Riemannian 
$n$-manifold $(M,g)$, 
$$B_0 \geq
\max
\left( {Vol(M)}^{-\frac{2}{n}},
\frac{{\mid
\mathcal{B}
\mid}^{-\frac{2}{n}}}{6n} \left( \frac{2}{n+2}+\frac{n-2}{{\lambda}_1} \right)
{\left(\frac{n+2}{2} \right)}^{\frac{2}{n}} \max_{x \in M}{S_g(x)} \right)$$
where
$S_g(x)$ is the scalar curvature of $g$ at $x$.
We now say that $u \in H_1^2(M)$, $u \not\equiv  0$ is an extremal function for the sharp $L^2$-
inequality $N(A_0(n),B_0)$ if 
$$(\int_{M}u^2dv_g)^{1+\frac{2}{n}}
=(A_0(n)\int_{M}{ \mid \nabla u \mid}_g^2
dv_{g} + B_0 \int_{M}u^2dv_g)(\int_{M}\mid u \mid
dv_g)^{\frac{4}{n}}$$
Such a study was carried out for sharp Sobolev inequalities by Djadli
and Druet in the very nice reference \cite{dd}. Though they are close
in their statement, these two questions, to know whether or not 
there exist extremal functions for sharp Sobolev inequalities and for
the sharp $L^2$-Nash inequality, are however distinct in nature. In 
consequence,
the problems we have to face here are very different from the one that appears
in \cite{dd}. The main
result of this article is the following :

\begin{theorem}
Let $(M,g)$ be a smooth compact Riemannian $n$-manifold. Let also $B_0$ be
as above. 
Assume that :
$$B_0>\frac{{\mid
\mathcal{B}
\mid}^{-\frac{2}{n}}}{6n} \left( \frac{2}{n+2}+\frac{n-2}{{\lambda}_1} \right)
{\left(\frac{n+2}{2} \right)}^{\frac{2}{n}} \max_{x \in M}{S_g(x)}$$
Then, there exist extremal functions of class $C^{1,a}(M)$ ( $0<a<1$ )
for the sharp $L^2$-Nash inequality.
\end{theorem}
We present here the main ideas of the proof of this theorem which is based 
on a precise study of a phenomenom of
concentration. Namely, for $B < B_0$, 
we prove the existence of an extremal function $u_B$
for inequality $N(A_B,B)$  where
 $$A_B = \inf \{A | \hbox{
  s.t. }  N(A,B) \hbox{ is true } \}> A_0(n)$$
We then let $B \to B_0$. 
Standard theory shows that there
exists $u \in H_1^2(M)$ such that $u_B \to u$ weakly in $H_1^2(M)$
when $B \to B_0$. We have to consider two cases. First, if $u \not\equiv 0
$, it is not difficult to prove that $u$ is
an extremal function for $N(A_0(n),B_0)$. If $u \equiv 0$, we prove that $u$
concentrates around a point $x$ of $M$. In other words,
$u_B \to 0$ when $B \to B_0$ in $C^0_{loc}(M-\{x\} )$ and for all $\delta >0$,
$$\lim_{B \to B_0} \frac{\int_{B(x,\delta)} u_B^2 dv_g}{\int_{M}
  u_B^2 dv_g}=1$$
Hence, if $\eta$ is a cut-off function
such that $\eta \equiv 1$ in a neighbourhood of $x$ and $\eta \equiv 0$ on
$M-B(x,\delta)$ where $\delta$ is small, $\eta u_B$ have almost the
same properties than $u_B$. 
Via exponential map at $x$, $\eta u_B$ can be seen
as a function on $\rR^n$ on which we have the standard optimal Nash inequality
$${\left(\int_{\rR^n} {(\eta u_B)}^2 dx\right)}^{1+\frac{2}{n}} \leq A_0(n) 
\int_{\rR^n} {\mid \nabla \eta u_B \mid }^2 dx 
{\left(\int_{\rR^n} {\mid \eta u_B \mid}dx\right)}^{\frac{4}{n}}$$
With the use of  Cartan's expansion of the metric around $x$ and precise
estimations of the concentration of $u_B$, 
these integrals can
be 
compared to the corresponding integrals on $(M,g)$. We get that 
$$\int_M {(\eta u_B)}^N dv_g \leq \alpha_B$$ 
where $\alpha_B$ is an expression involving integrals of $u_B$. 
Thanks to the Euler equation of $u_B$, we get that
$$\alpha_B' \leq \int_M {(\eta u_B)}^N dv_g $$
where $\alpha_B'$ is another  expression involving integrals of $u_B$.
The inegality $\alpha_B' \leq \alpha_B$ leads to 
$$B_0 \leq \frac{{\mid
\mathcal{B}
\mid}^{-\frac{2}{n}}}{6n} \left( \frac{2}{n+2}+\frac{n-2}{{\lambda}_1} \right)
{\left(\frac{n+2}{2} \right)}^{\frac{2}{n}} \max_{x \in M}{S_g(x)}$$
This gives the theorem.\\

 As a consequence of theorem 1, we immediately have :
\begin{cor}
Let $(M,g)$ be a smooth compact Riemannian $n$-manifold. We assume that 
$${Vol(M)}^{-\frac{2}{n}}>
\frac{{\mid
\mathcal{B}
\mid}^{-\frac{2}{n}}}{6n} \left( \frac{2}{n+2}+\frac{n-2}{{\lambda}_1} \right)
{\left(\frac{n+2}{2} \right)}^{\frac{2}{n}} \max_{x \in M}{S_g(x)}$$
Then, there exist extremal functions of class $C^{1,a}(M)$ ( $0<a<1$ )
for the sharp $L^2$-Nash inequality. In particular, this is the case if the
scalar curvature is nonpositive.
\end{cor}
For $n \geq 2$, 
the results obtained in \cite{h} on the existence of
extremal functions for the sharp $L^2$-Nash inequality are a consequence of
theorem 1. For $n=1$, we proved in \cite{h} that constant functions
are extremal functions for the sharp $L^2$-Nash inequality. At the moment,
we are not able to give examples manifolds such that there does not
exist extremal functions for the sharp $L^2$-Nash inequality. Hebey and
Vaugon prove in \cite{hv1} the existence of such manifolds in the case of 
Sobolev inequality. However, their proof strongly uses the conformal
invariance of their inequality and we do not know yet some other methods to
obtain this type of results.

%%%%%%%%%%%%%%%%%%%%%%%%%%%%%%%%%%%%%%%%%%%%%%%%%%%%%%%%%%%%%%%%%%%%%%%%%%

\section{Proof of theorem 1} 
Let $A_0(n)$ and $B_0$ be as in 
introduction. 
We  define  ${\alpha}_0=B_0 {A_0(n)}^{-1}$. For $\alpha >0$, we let also 
\begin{eqnarray*}
 I_{\alpha}(u) &  = & \frac{(\int_M {\mid \nabla u \mid}_g^2dv_g +
  ({\alpha}_0-{\alpha}) 
\int_M u ^2 dv_g )(\int_{M}{\mid u \mid}^{1+{\epsilon}_{\alpha}}
dv_g)^{\frac{4}{n(1+{\epsilon}_{\alpha})}}}
{{(\int_M u^2 dv_g)}^{1+\frac2n}} \\
 \Lambda & = & \{ u \in C^{\infty}(M) \hbox{ s.t. } \int_M u^2 dv_g = 1 \} 
\end{eqnarray*}
and  
$${\mu}_{\alpha}=\inf_{u \in \Lambda} I_{\alpha}(u)$$
where $\epsilon_{\alpha}$ is chosen such that 
\begin{eqnarray} \label{NR1}
\lim_{\alpha \to 0} \epsilon_{\alpha}=0, \hbox{ }
{\mu}_{\alpha} <  {A_0(n)}^{-1}  
\hbox{ and, } 
\lim_{\alpha \to 0} {\mu}_{\alpha}={A_0(n)}^{-1} 
\end{eqnarray}
Clearly there exists $u_{\alpha} \in H_1^2(M)$, $u_{\alpha} \geq 0$,
such that 
$$\int_M u_{\alpha}^2 dv_g = 1 \hbox{ and } \mu_{\alpha}=I_{\alpha}
(u_{\alpha})$$

\noindent We write now the Euler equation of
$u_{\alpha}$ to get that, in the sense of distributions :
$$2 A_{\alpha} {\Delta}_{g} u_{\alpha}+\frac{4}{n}
B_{\alpha}u_{\alpha}^{{\epsilon}_{\alpha}}= k_{\alpha}u_{\alpha}
\eqno{ (E_{\alpha})}$$
where ${\Delta}_{g}$ stands for the Laplacian with the minus sign convention
and :
\begin{eqnarray*}
  A_{\alpha} & = & {\left(
      \int_M{u_{\alpha}}^{1+{\epsilon}_{\alpha}}dv_g
\right)}^{\frac{4}{n(1+{\epsilon}_{\alpha})}} \\
  B_{\alpha} & = & \left(\int_M{ {\mid \nabla u_{\alpha} \mid}_{g}}^2 dv_g
  +  (\alpha_0-\alpha) \right)
{\left( \int_M{u_{\alpha}}^{1+{\epsilon}_{\alpha}}dv_g \right)}
^{\frac{4}{n(1+{\epsilon}_{\alpha})}-1} \\
  k_{\alpha} & = & \frac{4}{n} {\mu}_{\alpha} + 2\int_M{ \mid \nabla u_{\alpha}
  \mid}_g^2 dv_g {\left(\int_M{u_{\alpha}}^{1+{\epsilon}_{\alpha}}dv_g
\right)}^{\frac{4}{n(1+{\epsilon}_{\alpha})}} 
\end{eqnarray*}

\noindent By the Sobolev embedding theorem, $u_{\alpha} \in L^{\frac{2n}{n-2}}(M)$ and then, by
classical methods, $u_{\alpha} \in C^2(M)$. To prove the
theorem, we assume that there does not  exists extremal functions
for the sharp $L^2$-Nash inequality and show that 
$$B_0 \leq  \frac{{\mid
\mathcal{B}
\mid}^{-\frac{2}{n}}}{6n} \left( \frac{2}{n+2}+\frac{n-2}{{\lambda}_1} \right)
{\left(\frac{n+2}{2} \right)}^{\frac{2}{n}} \max_{x \in M}{S_g(x)}$$
As easily seen, 
the existence of extremal functions follows from an assumption like :
$$\liminf_{\alpha \to 0} \int_M {u_{\alpha}}^{1+ {\epsilon}_{\alpha}}
dv_g>0$$
Note that such an assumption implies that :
$$\int_M{ {\mid \nabla u_{\alpha} \mid}_{g}}^2 dv_g \leq C$$
In the following, we then assume that 
$$\lim_{\alpha \to 0} 
\int_M {u_{\alpha}}^{1+ {\epsilon}_{\alpha}} dv_g = 0$$
or, equivalently :
\begin{eqnarray} \label{r4}
\lim_{\alpha \to 0}  A_{\alpha} =0
\end{eqnarray}
\noindent Now, using $N(A_0(n), B_0)(u_{\alpha})$, we have :
$$\liminf_{\alpha \to 0} \int_M{ {\mid \nabla u_{\alpha} \mid}_{g}}^2 dv_g 
(\int_M{u_{\alpha}}^{1+{\epsilon}_{\alpha}}dv_g)^
{\frac{4}{n(1+{\epsilon}_{\alpha})}} \geq {A_0(n)}^{-1}$$
In addition, since ${\mu}_{\alpha}<{A_0(n)}^{-1}$, it is clear that : 
$$\limsup_{\alpha \to 0} \int_M{ {\mid \nabla u_{\alpha} \mid}_{g}}^2 dv_g 
(\int_M{u_{\alpha}}^{1+{\epsilon}_{\alpha}}dv_g)^
{\frac{4}{n(1+{\epsilon}_{\alpha})}} \leq A_0(n)^{-1}$$
As a consequence, one easily checks that :

\begin{eqnarray} \label{r1}
\lim_{\alpha \to 0} A_{\alpha} 
\int_M{ {\mid \nabla u_{\alpha} \mid}_{g}}^2 dv_g ={A_0(n)}^{-1}
\end{eqnarray}
\begin{eqnarray} \label{r2}
\lim_{\alpha \to 0} 
B_{\alpha} \int_M {u_{\alpha}}^{1+ {\epsilon}_{\alpha}} dv_g ={A_0(n)}^{-1}
\end{eqnarray}
\begin{eqnarray} \label{r3}
\lim_{\alpha \to 0} k_{\alpha}=(2+\frac4n){A_0(n)}^{-1}
\end{eqnarray}

\noindent The proof of the theorem proceeds in several  steps. 
Step 1 to  4 
are somehow  similar than what was done in \cite{h}. Note however that the
limits are not anymore limits as $\alpha \to \infty$. Step 5 is a
preparation
to the concluding step, step 6.

We let $a_{\alpha}={A_{\alpha}}^{\frac12}$. 
We let also $x_{\alpha}$ be a point of $M$ such that
 $u_{\alpha}(x_{\alpha})= {\parallel u_{\alpha} \parallel
 }_{\infty}$. In the following, $B(p,r)$ denotes the ball of center
 $p$ and radius $r$ in $\rR^n$ and $B_p(r)$ denotes the ball of center
 $p$ and radius $r$ in $M$. We assume in addition that bounded
 sequences are convergent, with no mention to the extracting of a
 subsequence, and write $C$ for positive constants that do not depend
 on $\alpha$. 

\begin{step}
For all $\delta >0$ : $\liminf_{\alpha \to 0} \frac{\int_{B_{x_{\alpha}}(\delta a_{\alpha})}
{u_{\alpha}}^{1+{\epsilon}_{\alpha}} dv_g }{\int_M {u_{\alpha}}^{1+{\epsilon}_{\alpha}}
dv_g }>0$
\end{step}

\noindent Let, for $x \in B(0, \delta) \subset \rR ^n$  : 
\begin{eqnarray*}
g_{\alpha}(x) & = & {(exp_{x_{\alpha}})}^*g(a_{\alpha}x) \\
{\varphi}_{\alpha}(x) & = & {\parallel u_{\alpha} \parallel  }_{\infty}^{-1}
 u_{\alpha}({exp}_{x_{\alpha}}(a_{\alpha}x))
\end{eqnarray*}

\noindent We easily get  :
$$ {\Delta}_{g_{\alpha}} {\varphi}_{\alpha}+\frac{2}{n}
{\parallel u_{\alpha} \parallel  }_{\infty}^{-1+{\epsilon}_{\alpha}}
B_{\alpha}{\varphi}_{\alpha}^{{\epsilon}_{\alpha}}
= \frac{k_{\alpha}}{2} {\varphi}_{\alpha}
\eqno{({\tilde E}_{\alpha})}$$
Since ${\Delta}_{g} u_{\alpha}(x_{\alpha}) \geq 0$, we get from $(E_{\alpha})$ and
(\ref{r3}) :
\begin{eqnarray} \label{r6}
{\parallel u_{\alpha} \parallel  }_{\infty}^{{\epsilon}_{\alpha}}B_{\alpha}
\leq C {\parallel u_{\alpha} \parallel  }_{\infty}
\end{eqnarray}
and since ${\parallel {\varphi}_{\alpha}
\parallel}_{L^{\infty}(B(0,\delta))} \leq 1$, we get from $({\tilde E}_{\alpha})$ :
$${\parallel {\Delta}_{g_{\alpha}} {\varphi}_{\alpha}
\parallel}_{L^{\infty}(B(0,\delta))} \leq C$$

\noindent By classical methods, it follows that, for $a \in ]0,1[$ : 
${\parallel {\varphi}_{\alpha}
\parallel}_{C^{1,a}B(0,\delta)}
\leq C$.
Hence, ${({\varphi}_{\alpha})}_{\alpha}$ is equicontinuous and by Ascoli's
theorem, there exists $\varphi$$\in$$C^{0}(B(0, \delta))$ such that 
${\varphi}_{\alpha} \to \varphi$ in $C^{0}(B(0, \delta))$ as $\alpha
\to 0$. We have :
\begin{eqnarray} \label{r8}
\varphi(0)= \lim_{\alpha \to 0} {\varphi}_{\alpha}(0) = 1
\end{eqnarray}

\noindent and also :

$$\int_{B(0, \delta)} {{\varphi}_{\alpha}}^{1+{\epsilon}_{\alpha}} dv_{g_{\alpha}}=
{\parallel u_{\alpha} \parallel  }_{\infty}^{-(1+{\epsilon}_{\alpha})}
{A_{\alpha}}^{-\frac{n}{2}} 
\int_{B_{x_{\alpha}}(\delta a_{\alpha})} {u_{\alpha}}^{1+{\epsilon}_{\alpha}} dv_g$$

$$={\parallel u_{\alpha} \parallel  }_{\infty}^{-(1+{\epsilon}_{\alpha})}
{A_{\alpha}}^{-\frac{n}{4}(1-{\epsilon}_{\alpha})} 
\frac{\int_{B_{x_{\alpha}}(\delta a_{\alpha})} {u_{\alpha}}^{1+{\epsilon}_{\alpha}}
dv_g}
{\int_M {u_{\alpha}}^{1+{\epsilon}_{\alpha}} dv_g}$$

\begin{eqnarray} \label{r9}
\leq
{\parallel u_{\alpha} \parallel  }_{\infty}^{-1}
{A_{\alpha}}^{-\frac{n}{4}} 
\frac{\int_{B_{x_{\alpha}}(\delta a_{\alpha})} {u_{\alpha}}^{1+{\epsilon}_{\alpha}}
dv_g}
{\int_M {u_{\alpha}}^{1+{\epsilon}_{\alpha}} dv_g}
\end{eqnarray}

\noindent Since ${\parallel u_{\alpha}
\parallel}_{\infty}^{{\epsilon}_{\alpha}} \geq 1$, (\ref{r6}) implies : ${\parallel u_{\alpha} \parallel 
}_{\infty} 
\geq C.  B_{\alpha}$
and since $A_{\alpha} \to 0$ as $\alpha \to 0$, (\ref{r2}) implies
that 
 $B_{\alpha} \geq C .
{A_{\alpha}}^{-\frac{n}{4}(1+ {\epsilon}_{\alpha})}
\geq C . {A_{\alpha}}^{-\frac{n}{4}}$. Inequality (\ref{r9}) then becomes :
$$\int_{B(0, \delta)} {{\varphi}_{\alpha}}^{1+ {\epsilon}_{\alpha}} dv_{g_{\alpha}}
\leq 
C 
\frac{\int_{B_{x_{\alpha}}(\delta a_{\alpha})}
{u_{\alpha}}^{1+{\epsilon}_{\alpha}} dv_g}
{\int_M {u_{\alpha}}^{1+{\epsilon}_{\alpha}} dv_g}$$
Moreover,
\begin{eqnarray} \label{r10} 
\int_{B(0, \delta)} {{\varphi}_{\alpha}}^{1+{\epsilon}_{\alpha}}
dv_{g_{\alpha}} \to
C >0
\end{eqnarray}
by (\ref{r8}) and since $g_{\alpha} \to \xi \hbox{ in } C^{1}(B)$
for every ball $B$ in $\rR^n$. Finally, we get :
$$\frac{
\int_{B_{x_{\alpha}}(\delta a_{\alpha})} {u_{\alpha}}^{1+{\epsilon}_{\alpha}} dv_g}
{\int_M {u_{\alpha}}^{1+{\epsilon}_{\alpha}} dv_g}
\geq
C >0$$
This ends the proof of step 1. Note that 
coming back to  (\ref{r9}) and (\ref{r10}), one easily gets that :
\begin{eqnarray} \label{r11}
\lim_{\alpha \to 0} {A_{\alpha}}^{\frac{n}{4}} 
 {\parallel u_{\alpha} \parallel }_{\infty}  =  C>0
\end{eqnarray}

\begin{step}
We recall that 
$$a_{\alpha}=A_{\alpha}^{\frac{1}{2}}={\left(
\int_{M} u_{\alpha}^{1+\epsilon_{\alpha}} dv_g 
\right)}^{\frac{2}{n(1+\epsilon_{\alpha})}}$$
Let ${(c_{\alpha})}_{\alpha}$ be a sequence of positive numbers such that :
$\frac{a_{\alpha}}{c_{\alpha}} \to 0$ as $\alpha \to 0$. Then :

$$\lim_{\alpha \to 0} \frac{\int_{B_{x_{\alpha}}(c_{\alpha})}
{u_{\alpha}}^{1+{\epsilon}_{\alpha}} dv_g}
{\int_M {u_{\alpha}}^{1+{\epsilon}_{\alpha}} dv_g}=1$$
\end{step}

\noindent Let $ \eta \in {C}^{\infty} (\rR) $ be such that :
\begin{eqnarray*}
(i) & \eta ([0, \frac{1}{2}])=\{ 1 \} \\
(ii) & \eta ([1, +\infty[)=\{ 0 \} \\
(iii) & 0 \leq \eta \leq 1
\end{eqnarray*}

\noindent For $ k \in \nN $, we let : ${\eta}_{\alpha, k}(x)={\left( \eta (
{{c}_{\alpha}}^{-1}  {d}_{g}(x, {x}_{\alpha})) \right)}^{{2}^{k}}$ where 
$d_g$ denotes the distance for $g$.
Multiplying $({E}_{\alpha})$ by ${{\eta}_{\alpha, k}}^{2}
{u}_{\alpha}$
 and integrating
over $M$ gives :

\begin{eqnarray*}
2 {A}_{\alpha} \int_{M} { \mid \nabla {\eta}_{\alpha, k} 
{u}_{\alpha} \mid }_g^{2} {dv}_{g}
- 2 {A}_{\alpha} \int_{M} 
{\mid \nabla {\eta}_{\alpha, k}  \mid}_g^{2} {{u}_{\alpha}}^{2} {dv}_{g} +
\frac{4}{n} {B}_{\alpha} \int_{M} {{\eta}_{\alpha, k}}^{2}
{{u}_{\alpha}}^{1+{\epsilon}_{\alpha}} dv_{g}
\end{eqnarray*}
\begin{eqnarray} \label{r13}
= {k}_{\alpha} \int_{M} {({\eta}_{\alpha, k} {u}_{\alpha})}^{2} {dv}_{g}
\end{eqnarray}
Using $N(A_0(n)+\epsilon, B_{\epsilon})({\eta}_{\alpha, k} u_{\alpha})$, one easily
checks :

\begin{eqnarray*}
2A_{\alpha} \int_{M} {\mid \nabla {\eta}_{\alpha, k} u_{\alpha} \mid}_g^2 dv_g -
2A_{\alpha} \int_{M} {\mid \nabla {\eta}_{\alpha, k}  \mid}_g^2
{u_{\alpha}}^2 dv_g +
\frac{4}{n} B_{\alpha} \int_M {{\eta}_{\alpha, k}}^2
{u_{\alpha}}^{1+{\epsilon}_{\alpha}} dv_g 
\end{eqnarray*}
\begin{eqnarray} \label{r14}
\nonumber \leq
k_{\alpha} \Big( (A_0(n)+\epsilon) 
\int_{M} {\mid \nabla {\eta}_{\alpha, k} u_{\alpha} \mid}_g^{2} dv_{g} 
{(\int_{M} {({\eta}_{\alpha, k} u_{\alpha})}^{1+{\epsilon}_{\alpha}}
dv_{g})}^{\frac{4}{n(1+{\epsilon}_{\alpha})}} +\\
{B}_{\epsilon} \int_{M} {({\eta}_{\alpha, k} u_{\alpha})}^2 dv_g 
{(\int_{M} {({\eta}_{\alpha, k} u_{\alpha})}^{1+{\epsilon}_{\alpha}}
dv_g)}^{\frac{4}{n(1+{\epsilon}_{\alpha})}} \Big)^{\frac{n}{n+2}}
\end{eqnarray}

\noindent Moreover, with the assumption on $(c_{\alpha})_{\alpha}$ : 

$${\mid \nabla {\eta}_{\alpha, k}  \mid}_g^2 \leq 
\frac{C}{{c_{\alpha}}^2}
\Rightarrow
\lim_{\alpha \to 0} A_{\alpha} \int_{M} 
{\mid \nabla {\eta}_{\alpha, k}  \mid}_g^2 {u_{\alpha}}^2 dv_g=0$$
Now, let :

\begin{eqnarray*}
 {\lambda}_{k} & = & \lim_{\alpha \to 0} \frac {\int_M {{\eta}_{\alpha, k}}^2
{u_{\alpha}}^{1+{\epsilon}_{\alpha}} dv_g }
{\int_{M} {u_{\alpha}}^{1+{\epsilon}_{\alpha}} dv_g } \\                  
{\tilde {\lambda}}_k & = & \lim_{\alpha \to 0} \frac {\int_{M} {({\eta}_{\alpha, k}
u_{\alpha})}^{1+{\epsilon}_{\alpha}} dv_g }
{\int_{M} {u_{\alpha}}^{1+{\epsilon}_{\alpha}} dv_g }
\end{eqnarray*}

\noindent From the definition of ${\eta}_{\alpha, k}$, we get, for all $k \in \nN$ :

\begin{eqnarray} \label{r15}
{\lambda}_{k+1} \leq {\tilde {\lambda}}_{k+1}
\leq {\lambda}_{k} \leq {\tilde{\lambda}}_{k}
\leq \mu=\lim_{\alpha \to 0} \frac{\int_{B_{x_{\alpha}}(c_{\alpha})}  
{u_{\alpha}}^{1+{\epsilon}_{\alpha}} dv_g}{\int_M {u_{\alpha}}^{1+{\epsilon}_{\alpha}}
dv_g }
\end{eqnarray}

\noindent and, by step 1 : 
\begin{eqnarray} \label{r16}
\exists C>0
\hbox{ s.t.   }
\forall k
\in
\nN, {\lambda}_{k}
\geq C
\end{eqnarray}

\noindent Let us now prove that :
${\lambda}_{k} \leq {{\tilde {\lambda}}_{k}}^2$.
Let $L_k=\lim_{\alpha \to 0} A_{\alpha}
\int_{M} {\mid \nabla {\eta}_{\alpha, k} u_{\alpha} \mid}_g^2 dv_g $.
Note that (\ref{r2}) and (\ref{r3}) imply :

\begin{eqnarray*}
\lim_{\alpha \to 0} B_{\alpha} \int_M {{\eta}_{\alpha, k}}^2
{u_{\alpha}}^{1+{\epsilon}_{\alpha}} dv_g = {\lambda}_k {A_0(n)}^{-1} \\
\end{eqnarray*}
and 
\begin{eqnarray*}
k_{\alpha} \int_{M} {({\eta}_{\alpha, k} u_{\alpha})}^2 dv_g \leq
C
\end{eqnarray*}

\noindent In particular,  (\ref{r13}) gives : $L_{k} < + \infty$.
We also clearly have  by (\ref{r1}) and (\ref{r2}) :
$$\lim_{\alpha \to 0} \int_{M} {\mid \nabla {\eta}_{\alpha, k} u_{\alpha} \mid}_g^2 dv_g 
{(\int_{M} {({\eta}_{\alpha, k} u_{\alpha})}^{1+{\epsilon}_{\alpha}} dv_g
)}^{\frac{4}{n(1+{\epsilon}_{\alpha})}} = L_k {\tilde{\lambda}}_{k}^{\frac{4}{n}}$$

\noindent Equation (\ref{r14}) then leads to :
$$2 L_k+ \frac{4}{n} {A_0(n)}^{-1} {\lambda}_{k}
\leq (2+ \frac{4}{n}) {A_0(n)}^{-1} 
{((A_0(n) + \epsilon) L_k {{\tilde{\lambda}}_{k}}^{\frac{4}{n}})}^{\frac{n}{n+2}}$$

\noindent If ${{\tilde L}_k=A_0(n)} L_k$, we obtain, since $\epsilon$ was arbitrary :

\begin{eqnarray*}
2{\tilde L}_k+\frac{4}{n} {\lambda}_{k} \leq (2+ \frac{4}{n}) {\tilde
L}_{k}^{\frac{n}{n+2}}
{{\tilde {\lambda}}_{k}}^{\frac{4}{n+2}}
\end{eqnarray*}

\noindent Let now, for $x, y, z$ :
$f(x, y, z)=(2+ \frac{4}{n}) x^{\frac{n}{n+2}} y^{\frac{4}{n+2}}
- (\frac{4}{n}z+2x)$. Differentiating in $x$, we see that
$\forall x, y, z>0, f(x, y, z) \leq f(y^2, y, z)$,
and then : $f({\tilde L}_k, {\tilde{\lambda}}_k, {\lambda}_k) \leq 
f({{\tilde{\lambda}}_k}^2, {\tilde{\lambda}}_k, {\lambda}_k) 
= \frac{4}{n}({{\tilde{\lambda}}_k}^2-{\lambda}_k)$.
We then get : ${\lambda}_{k} \leq {{\tilde {\lambda}}_{k}}^2$. 
Now, from (\ref{r15}), (\ref{r16}), we get : $ \forall N \in \nN ,
\hbox{ }  0<C \leq
{{\lambda}_{0}}^{N} \leq \mu$.
Since $\mu \leq 1$, we have $\mu=1$ which proves
step 2. Note that 
we have also proved that $\tilde{L}_k=1$ for all $k$. As one can check,
we have then :
\begin{eqnarray} \label{az3}
\lim_{\alpha \to 0} \frac{\int_{B_{x_{\alpha}}(c_{\alpha})}
{\mid \nabla u_{\alpha} \mid}_g^2 dv_g}{\int_M
{\mid \nabla u_{\alpha} \mid}_g^2 dv_g}=1
\end{eqnarray}
As a consequence,  we easily get from (\ref{r13}) :
\begin{eqnarray} \label{az2}
\lim_{\alpha \to 0} \int_{B_{x_{\alpha}}(c_{\alpha})}
u_{\alpha}^2 dv_g=\lim_{\alpha \to 0} \frac{\int_{B_{x_{\alpha}}(c_{\alpha})}
u_{\alpha}^2 dv_g}{\int_M u_{\alpha}^2 dv_g}=1
\end{eqnarray}

\begin{step}
There exists $C>0$ such that, for all $x \in M$ : 
$$u_{\alpha}(x){d(x, x_{\alpha})}^{\frac{n}{2}} \leq C$$
where $d$ denotes the distance for $g$.
\end{step}

\noindent We proceed by contradiction. 
We suppose that the following assumption is true :
$$\exists y_{\alpha} \in M \hbox{ s.t.   }
\lim_{\alpha \to 0} 
u_{\alpha}(y_{\alpha}){d(y_{\alpha}, x_{\alpha})}^{\frac{n}{2}}=
+\infty \eqno{(H)}$$
Let :
$$v_{\alpha}=u_{\alpha}(y_{\alpha}){d(y_{\alpha}, x_{\alpha})}^{\frac{n}{2}}$$
We can assume that :
$$v_{\alpha} = {\parallel u_{\alpha}(.){d(., x_{\alpha})}^{\frac{n}{2}} 
\parallel}_{\infty}$$

\noindent First, we prove that, if $\nu$ is small enough :

\begin{eqnarray} \label{r19}
B_{y_{\alpha}}({u_{\alpha}(y_{\alpha})}^{-\frac{2}{n}}) \cap
B_{x_{\alpha}}(a_{\alpha}{v_{\alpha}}^{\nu})= \emptyset
\end{eqnarray}
It is here sufficient to show that $d(x_{\alpha}, y_{\alpha}) \geq
{u_{\alpha}(y_{\alpha})}^{-\frac{2}{n}}+a_{\alpha}{v_{\alpha}}^{\nu}$,
or, equivalently that ${v_{\alpha}}^{\frac{2}{n}-\nu} \geq {v_{\alpha}}^{-\nu}
+a_{\alpha}{u_{\alpha}(y_{\alpha})}^{\frac{2}{n}}$.
If $\nu < \frac{2}{n}$, from $(H)$, we get that 
${v_{\alpha}}^{\frac{2}{n}-\nu} 
\to  +\infty$ and $ {v_{\alpha}}^{-\nu} \to 0$ as $\alpha \to 0$.  
Hence, it still has to be proved that 
$a_{\alpha}{u_{\alpha}(y_{\alpha})}^{\frac{2}{n}} 
\leq C$. 
We have $a_{\alpha}{u_{\alpha}(y_{\alpha})}^{\frac{2}{n}} \leq  a_{\alpha}
{\parallel u_{\alpha} \parallel}_{\infty}^{\frac{2}{n}}$.
Since ${a}_{\alpha}={{A}_{\alpha}}^{\frac{1}{2}}$ 
and by (\ref{r11}), this gives : $a_{\alpha}
{\parallel u_{\alpha} \parallel}_{\infty}^{\frac{2}{n}} \leq C$ .
Equation (\ref{r19}) then follows.  
We let now, for $x \in B(0,1)$ :

\begin{eqnarray*}
h_{\alpha}(x) & = & {(exp_{y_{\alpha}})}^*g(l_{\alpha}x) \\
{\psi}_ {\alpha}(x) & = & {u_{\alpha}(y_{\alpha})}^{-1}
 u_{\alpha}({exp}_{y_{\alpha}}(l_{\alpha}x))
\end{eqnarray*}
where :
$$l_{\alpha}= {\parallel u_{\alpha} \parallel}_{\infty}^{-\frac{n+4}{2n}}
{u_{\alpha}(y_{\alpha})}^{\frac{1}{2}}$$
On $B(0, 1)$, we have :
$${\Delta}_{h_{\alpha}} {\psi}_ {\alpha}=
\frac{k_{\alpha} {\parallel u_{\alpha} \parallel}_{\infty}^{-(1+\frac{4}{n})}
u_{\alpha}(y_{\alpha})}{2 A_{\alpha}} 
{\psi}_ {\alpha} -
\frac{2 B_{\alpha} {\parallel u_{\alpha} \parallel}_{\infty}^{-(1+\frac{4}{n})}
{u_{\alpha}(y_{\alpha})}^{{\epsilon}_{\alpha}}}{n A_{\alpha}}
{{\psi}_ {\alpha}}^{{\epsilon}_{\alpha}} \eqno{({E_{\alpha}}')}$$
Moreover :
\begin{eqnarray} \label{r21}
h_{\alpha} \to \xi \hbox{ in } C^1(B(0, 1)) \hbox{ as } \alpha \to 0
\end{eqnarray}

\noindent We have ${\parallel u_{\alpha} \parallel}_{L^{\infty}
(B_{y_{\alpha}}(u_{\alpha}(y_{\alpha})^{-\frac{2}{n}}))}
\leq C.  u_{\alpha}(y_{\alpha})$.
To see this, note that, by the definition of $y_{\alpha}$, we have for all
$x \in B_{y_{\alpha}}( {u_{\alpha}(y_{\alpha})}^{-\frac{2}{n}})$ :
\begin{eqnarray} \label{r22}
u_{\alpha}(y_{\alpha}){d(x_{\alpha},y_{\alpha})}^{\frac{n}{2}}
\geq u_{\alpha}(x){d(x_{\alpha},x)}^{\frac{n}{2}}
\end{eqnarray}

\noindent Moreover, since $x \in B_{y_{\alpha}}(
{u_{\alpha}(y_{\alpha})}^{-\frac{2}{n}})$ :
$$d(y_{\alpha},x) \leq 
{u_{\alpha}(y_{\alpha})}^{-\frac{2}{n}}$$
and, by $(H)$ : ${u_{\alpha}(y_{\alpha})}^{-\frac{2}{n}} \leq \frac{1}{2}
d(x_{\alpha},y_{\alpha})$.
So we have :
$$d(x,x_{\alpha}) \geq d(x_{\alpha},y_{\alpha})-d(x,y_{\alpha}) 
\geq d(x_{\alpha},y_{\alpha})-{u_{\alpha}(y_{\alpha})}^{-\frac{2}{n}} \geq
\frac{1}{2} d(x_{\alpha},y_{\alpha})$$
 
\noindent Coming back to (\ref{r22}), the result follows immediately.
Since $l_{\alpha} \leq {u_{\alpha}(y_{\alpha})}^{-\frac{2}{n}}$, we
then have 
${\parallel {\psi}_{\alpha} \parallel}_{L^{\infty}(B(0,1))} \leq C$. 
From (\ref{r6}), (\ref{r11}) and the fact that, by
(\ref{r2}), $B_{\alpha} {A_{\alpha}}^{\frac{n}{4}(1+{\epsilon}_{\alpha})} \to
C >0$ as $\alpha \to 0$, we get 
\begin{eqnarray} \label{nr3}
\lim_{\alpha \to 0} 
{\parallel u_{\alpha} \parallel}_{\infty}^{{\epsilon}_{\alpha}} = C
\end{eqnarray}

\noindent Now, from (\ref{r6}), (\ref{r11}) and (\ref{nr3}), we see that $(E_{\alpha}')$ has
bounded coefficients and then :
$${\parallel {\Delta}_{h_{\alpha}} {\psi}_{\alpha} \parallel}
_{L^{\infty}(B(0, 1))} \leq C$$

\noindent As in step 1, 
we get the existence of ${\psi} \in C^0(B(0, 1))$ such that,
up to a subsequence :
$${\psi}_{\alpha} \to {\psi} \hbox{ in } C^0(B(0, 1)) \hbox{ as }
\alpha \to 0$$
Here, ${\psi}$ is such that ${\psi}(0)=1$ and then :

\begin{eqnarray} \label{r23}
\int_{B(0, 1)} {\psi} dx >0
\end{eqnarray}

\noindent However, by (\ref{r21}) :
$$\int_{B(0, 1)} {\psi} dx=
\lim_{\alpha \to 0} \int_{B(0, 1)} {{\psi}_{\alpha}}^{1+{\epsilon}_{\alpha}} dv_{h_{\alpha}}$$

\noindent and, as one can check :
$$\int_{B(0, 1)} {{\psi}_{\alpha}}^{1+{\epsilon}_{\alpha}}
dv_{h_{\alpha}}=\beta_{\alpha}$$
where 
$$\beta_{\alpha}={A_{\alpha}}^{\frac{n}{4}(1+{\epsilon}_{\alpha})}
{u_{\alpha}(y_{\alpha})}^{-(1+{\epsilon}_{\alpha})} {l_{\alpha}}^{-n}
\left( \frac
{\int_{B_{y_{\alpha}}(l_{\alpha})} {u_{\alpha}}^{1+{\epsilon}_{\alpha}} dv_g}
{{A_{\alpha}}^{\frac{n}{4}(1+{\epsilon}_{\alpha})}} \right) $$

\noindent If we prove that $\lim_{\alpha \to 0} \beta_{\alpha} =0$, we get a
contradiction with (\ref{r23}) which ends the proof of step 3.
First, let 
$$m_{\alpha}= \frac{u_{\alpha}(y_{\alpha})}{{\parallel u_{\alpha}
    \parallel}_{\infty}}$$
Clearly, by (\ref{r11})  : 
$$\beta_{\alpha} \leq C m_{\alpha}^{-(\frac{n}{2}+1)}  
\left( \frac
{\int_{B_{y_{\alpha}}( u_{\alpha}(l_{\alpha}))} {u_{\alpha}}^{1+{\epsilon}_{\alpha}} dv_g}
{\int_M {u_{\alpha}}^{1+{\epsilon}_{\alpha}} dv_g} \right)$$
By step 2 and (\ref{r19}), 
\begin{eqnarray} \label{az1}
\lim_{\alpha \to 0} \left( \frac
{\int_{B_{y_{\alpha}}( u_{\alpha}(y_{\alpha})^{-\frac{2}{n}})} {u_{\alpha}}^{1+{\epsilon}_{\alpha}} dv_g}
{\int_M {u_{\alpha}}^{1+{\epsilon}_{\alpha}} dv_g} \right) = 0
\end{eqnarray}
If $m_{\alpha} \geq C >0$, we have $\beta_{\alpha} \to 
0$ as $\alpha \to 0$. 
Hence, we assume that $\lim_{\alpha \to 0} m_{\alpha} =0$.  
We now proceed  by induction to prove that :
$$\lim_{\alpha \to 0} m_{\alpha}^{ - {\left( \frac{n+3} {n+2} \right)}^{k}}  
\int_{  B_{y_{\alpha}}(2^{-k} u_{\alpha}(y_{\alpha})^{-\frac{2}{n}})}
 {u_{\alpha}}^2 dv_{g}= 0 \eqno{(H_{k})}$$
First, we prove that $(H_0)$ is true. We proved before that
$${\parallel u_{\alpha} \parallel}_{L^{\infty}
(B_{y_{\alpha}}(u_{\alpha}(y_{\alpha})^{-\frac{2}{n}}))}
\leq C.  u_{\alpha}(y_{\alpha})$$ 
Hence, we have, noting that $u_{\alpha}(y_{\alpha}) \to \infty$ as
$\alpha \to 0$ :
$$\int_{  B_{y_{\alpha}}( u_{\alpha}(y_{\alpha})^{-\frac{2}{n}})} 
u_{\alpha}^2 dv_g
\leq C u_{\alpha}(y_{\alpha}) 
\int_{  B_{y_{\alpha}}( u_{\alpha}(y_{\alpha})^{-\frac{2}{n}})} 
u_{\alpha}^{1+\epsilon_{\alpha}} dv_g$$
$$ \leq C m_{\alpha} 
{\parallel u_{\alpha} \parallel}_{\infty}
\int_{  B_{y_{\alpha}}( u_{\alpha}(y_{\alpha})^{-\frac{2}{n}})} 
u_{\alpha}^{1+\epsilon_{\alpha}} dv_g$$
By (\ref{r11}) and (\ref{az1})
$$\lim_{\alpha \to 0} {\parallel u_{\alpha} \parallel}_{\infty} \int_{  B_{y_{\alpha}}( u_{\alpha}(y_{\alpha})^{-\frac{2}{n}})} 
u_{\alpha}^{1+\epsilon_{\alpha}} dv_g=0$$
$(H_0)$ then follows.
Let now ${\epsilon}_k={\left( \frac{n+3}{n+2} \right)}^k$ and suppose that $(H_k)$
is true. 
Let us prove that $(H_{k+1})$ is true. 
Let ${\eta}_{\alpha, k} (x) = \eta (u_{\alpha}(y_{\alpha})^{\frac{2}{n}} 2^{k}
d_g(x, y_{\alpha}))$  where $\eta$ is defined as in step 2.  
Multiplying $(E_{\alpha})$ by 
$$\frac{u_{\alpha} {({\eta}_{\alpha, k})}^2} {m_{\alpha}^{{\epsilon}_k}}$$
and integrating over M, we obtain :

$$2 {A_{\alpha}} m_{\alpha}^{-\epsilon_k}
\int_{M} { \mid \nabla {\eta}_{\alpha, k} 
{u}_{\alpha} \mid }_g^{2} {dv}_{g}
-2 {A_{\alpha}} m_{\alpha}^{-\epsilon_k}
\int_M {\mid \nabla {\eta}_{\alpha, k}  \mid}_g^{2} 
{{u}_{\alpha}}^{2} {dv}_{g}$$ 
\begin{eqnarray} \label{nr8}
+\frac{4}{n} {B}_{\alpha} {m_{\alpha}}^{-{\epsilon}_k}   \int_{M}
{{\eta}_{\alpha, k}}^{2} {{u}_{\alpha}}^{1+{\epsilon}_{\alpha}} dv_{g}
= {k}_{\alpha} { {m_{\alpha}}^{-{\epsilon}_k} }
\int_{M} {({\eta}_{\alpha, k} {u}_{\alpha})}^{2} {dv}_{g}
\end{eqnarray}

\noindent By $(H_k)$ :

$$2 {A_{\alpha}} m_{\alpha}^{-\epsilon_k}
\int_M {\mid \nabla {\eta}_{\alpha, k}  \mid}_g^{2}
 {{u}_{\alpha}}^{2} {dv}_{g}$$ 
$$\leq C A_{\alpha} 
u_{\alpha}(y_{\alpha})^{\frac{4}{n}} m_{\alpha}^{-\epsilon_k} 
\int_{B_{y_{\alpha}}(2^{-k} u_{\alpha}(y_{\alpha})^{-\frac{2}{n}})}
u_{\alpha}^2 dv_g\leq C A_{\alpha}
u_{\alpha}(y_{\alpha})^{\frac{4}{n}}$$
Moreover, by (\ref{r11}), 
$ A_{\alpha}
u_{\alpha}(y_{\alpha})^{\frac{4}{n}} = A_{\alpha} m_{\alpha}^{\frac{4}{n}}
{\parallel u_{\alpha} \parallel}_{\infty}^{\frac{4}{n}} 
\leq C. m_{\alpha}^{\frac{4}{n}} \to 0$ as $\alpha \to 0$.
We have also, by $(H_k)$ and (\ref{r3}) :
$$\lim_{\alpha \to 0} {{k}_{\alpha}} { {m_{\alpha}}^{-{\epsilon}_k} }
\int_{M} {({\eta}_{\alpha, k} {u}_{\alpha})}^{2} {dv}_{g} = 0$$

\noindent Therefore, (\ref{nr8}) gives :

\begin{eqnarray} \label{o25}
2 {A_{\alpha}} 
\int_{M} { \mid \nabla {\eta}_{\alpha, k} 
{u}_{\alpha} \mid }_g^{2} {dv}_{g} \leq C.m_{\alpha}^{\epsilon_k} 
\end{eqnarray}
\begin{eqnarray*} 
\frac{4  }{n }{B}_{\alpha} \int_{M}
{{\eta}_{\alpha, k}}^{2} {{u}_{\alpha}}^{1+{\epsilon}_{\alpha}} dv_{g} 
 \leq C.m_{\alpha}^{\epsilon_k} 
\end{eqnarray*}
Up to replacing $\eta_{\alpha,k}$ by $\sqrt{{\eta}_{\alpha, k}}$, with
the same arguments, we also have :
\begin{eqnarray} \label{o26}
\frac{4  }{n }{B}_{\alpha} \int_{M}
{{\eta}_{\alpha, k}}^{1+\epsilon_{\alpha}} {{u}_{\alpha}}^{1+{\epsilon}_{\alpha}} dv_{g} 
 \leq C.m_{\alpha}^{\epsilon_k} 
\end{eqnarray}

\noindent Moreover, using $N(A, B)({\eta}_{\alpha, k} u_{\alpha} )$,
one easily checks that :

\begin{eqnarray*}
{ \left( \int_{M} {({\eta}_{\alpha, k} {u}_{\alpha})}^{2} {dv}_{g} 
 \right) }^{\frac{n+2}{n}}
\leq 
A . 
{ \int_{M} { \mid \nabla {\eta}_{\alpha, k} 
{u}_{\alpha} \mid }_g^{2} {dv}_{g} 
{ \left( \int_{M}
{({\eta}_{\alpha, k} {u}_{\alpha})}^{1+{\epsilon}_{\alpha}} dv_{g} \right) }   
^{\frac{4}{n(1+{\epsilon}_{\alpha})}} }
\end{eqnarray*}

\begin{eqnarray*}
+ B .  \int_{M} {({\eta}_{\alpha, k} {u}_{\alpha})}^{2} {dv}_{g}
{ \left( \int_{M}
{({\eta}_{\alpha, k} {u}_{\alpha})}^{1+{\epsilon}_{\alpha}} dv_{g} \right) }   
^{\frac{4}{n(1+{\epsilon}_{\alpha})}} 
\end{eqnarray*} 
Clearly,  we have in fact that : 
$${ \left( \int_{M} {({\eta}_{\alpha, k} {u}_{\alpha})}^{2} {dv}_{g} 
 \right) }^{\frac{n+2}{n}}
\leq 
C . 
{ \int_{M} { \mid \nabla {\eta}_{\alpha, k} 
{u}_{\alpha} \mid }_g^{2} {dv}_{g} 
{ \left( \int_{M}
{({\eta}_{\alpha, k} {u}_{\alpha})}^{1+{\epsilon}_{\alpha}} dv_{g} \right) }   
^{\frac{4}{n(1+{\epsilon}_{\alpha})}} }$$
$$\leq \frac{C}{A_{\alpha}
  B_{\alpha}^{\frac{4}{n(1+\epsilon_{\alpha})}}}
\left(  \int_{M} { \mid \nabla {\eta}_{\alpha, k} 
{u}_{\alpha} \mid }_g^{2} {dv}_{g} A_{\alpha} \right)
{\left( B_{\alpha} \int_{M}
{({\eta}_{\alpha, k} {u}_{\alpha})}^{1+{\epsilon}_{\alpha}} dv_{g} \right) }   
^{\frac{4}{n(1+{\epsilon}_{\alpha})}}$$
Using (\ref{o25}) and (\ref{o26}), we get 
$${ \left( \int_{M} {({\eta}_{\alpha, k} {u}_{\alpha})}^{2} {dv}_{g} 
 \right) }^{\frac{n+2}{n}} \leq \frac{C}{A_{\alpha}
  B_{\alpha}^{\frac{4}{n(1+\epsilon_{\alpha})}}}
.m_{\alpha}^{(1+\frac{4}{n(1+\epsilon_{\alpha})})\epsilon_k}$$
By (\ref{r2}), $A_{\alpha}
  B_{\alpha}^{\frac{4}{n(1+\epsilon_{\alpha})}} \geq C >0$. Since :
$$\int_{B_{y_{\alpha}}(2^{-(k+1)}
u_{\alpha}(y_{\alpha})^{-\frac{2}{n}})} u_{\alpha}^2 dv_g \leq 
\int_{M} {({\eta}_{\alpha, k} {u}_{\alpha})}^{2} {dv}_{g} $$
$(H_{k+1})$ then follows. As a consequence,
 $(H_k)$ is true for all $k$.  Coming back to (\ref{o26}) , we get
 that, for all $k$ :
$$\lim_{\alpha \to 0} m_{\alpha}^{-\epsilon_k} B_{\alpha}
 \int_{B_{y_{\alpha}}(2^{-k} u_{\alpha}(y_{\alpha})^{-\frac{2}{n}})}
 {{u}_{\alpha}}^{1+{\epsilon}_{\alpha}} dv_{g} 
 =0$$
Using the fact that $\lim_{\alpha \to 0}l_{\alpha}u_{\alpha}
  (y_{\alpha})^{\frac{2}{n}}=0$  and choosing $k$ such that 
$\epsilon_k \geq \frac{n}{2}+1$,
we get :  $\lim_{\alpha \to 0} \beta_{\alpha}=0$ which ends the proof
of step 3.

\begin{step} For all $c,k>0$,
we have :
\begin{eqnarray} \label{r31}
\lim_{\alpha \to 0} {A_{\alpha}}^{-k} \int_{M-B_{x_{\alpha}}( c)}
{u_{\alpha}}^2 dv_g
& = &  0
\end{eqnarray}
\begin{eqnarray} \label{r32}
\lim_{\alpha \to 0} {A_{\alpha}}^{-k} \int_{M-B_{x_{\alpha}}(c)} {\mid \nabla u_{\alpha} 
\mid}_g^2 dv_g & = & 0
\end{eqnarray}
\begin{eqnarray} \label{r33}
\lim_{\alpha \to 0} {A_{\alpha}}^{-k} 
\int_{M-B_{x_{\alpha}}(c)} {u_{\alpha}}^{1+{\epsilon}_{\alpha}} dv_g & = & 0
\end{eqnarray}
\end{step}

\noindent Let $r_{\alpha}(x)=d_g(x, x_{\alpha})$ 
and let $\delta \in ]0, \frac{n}{4}[$. Using
step 3, we have :
$${A_{\alpha}}^{-\delta} \int_{M-B_{x_{\alpha}}(c)}
{u_{\alpha}}^2 dv_g
\leq C.
{A_{\alpha}}^{-\delta}
\int_{M-B_{x_{\alpha}}(c)} {u_{\alpha}}^{1+{\epsilon}_{\alpha}} 
{r_{\alpha}}^{-\frac{n}{2}(1-{\epsilon}_{\alpha})} dv_g$$
$$\leq C. 
{A_{\alpha}}^{-\delta}
\int_{M-B_{x_{\alpha}}(c)} {u_{\alpha}}^{1+{\epsilon}_{\alpha}}  dv_g$$
Recall the definition of $A_{\alpha}$ to get :
$$\lim_{\alpha \to 0}  {A_{\alpha}}^{-\delta} \int_{M-B_{x_{\alpha}}(c)}
{u_{\alpha}}^2 dv_g = 0$$

\noindent Mimicking what we have done in the proof of step 3, we prove  by
induction that, for all k :
$$\lim_{\alpha \to 0} { A_{\alpha} }^{ - {\left( \frac{n+3} {n+2} \right)}^{k} \delta}  \int_{ M-
 B_{x_{\alpha}}(2^{k} c)} {u_{\alpha}}^2 dv_{g} = 0 $$
This  gives (\ref{r31}). Following the arguments used in the proof 
of step 3, one easily
gets (\ref{r32}) and (\ref{r33}) from (\ref{o25}) and (\ref{o26}).
Now,  we set, for $c>0$ small, 
${\eta}_{\alpha}=\eta (c^{-1} r_{\alpha})$ where $\eta$ is as
above. 
We also define :
\begin{eqnarray*}
r_{\nabla} & = & \frac{\int_M 
{\mid \nabla u_{\alpha} {\eta}_{\alpha} \mid}_g^2 
R_{ij}(x_{\alpha})x^{i}x^{j}dv_g}
{\int_M {\mid \nabla u_{\alpha} {\eta}_{\alpha} \mid}_g^2 dv_g}\\
r_1 & = & \frac{\int_M {( u_{\alpha} {\eta}_{\alpha} )}^
{1+{\epsilon}_{\alpha}} 
R_{ij}(x_{\alpha})x^{i}x^{j}dv_g}
{\int_M {( u_{\alpha} {\eta}_{\alpha} )}^{1+{\epsilon}_{\alpha}}  dv_g}\\
r_2 & = & \frac{\int_M {( u_{\alpha} {\eta}_{\alpha} )}^2
R_{ij}(x_{\alpha})x^{i}x^{j}dv_g}
{\int_M {( u_{\alpha} {\eta}_{\alpha} )}^2  dv_g}
\end{eqnarray*}
where $(x^1,..,x^n)$ are exponential coordinates.

\begin{step}
We have 
$$\lim_{\alpha \to 0} \frac{-\frac{1}{6} \left( -r_{\nabla}+(1+\frac{2}{n})r_2-
\frac{4}{n(1+{\epsilon}_{\alpha})} r_1 \right)}
{A_{\alpha}}$$
\begin{eqnarray} \label{s1}
=\frac{{\mid \mathcal{B} \mid}^{-\frac{2}{n}}}{6n}
\left( 
\frac{2}{n+2}+\frac{n-2}{\lambda}_1 \right) {\left( \frac{n+2}{2}
\right)}^{\frac{2}{n}} {S_g(x_0)} 
\end{eqnarray}
\end{step}

\noindent We come back to the notations of step 1. 
Let :
$$C_0=\lim_{\alpha \to 0} {\parallel u_{\alpha} \parallel}_{\infty}^{-1}
A_{\alpha}^{-\frac{n}{4}} \hbox{ and }
{\tilde{C}}_0=\lim_{\alpha \to 0} A_{\alpha}^{{\epsilon}_{\alpha}}$$

\noindent Note that,   by (\ref{r11}) and (\ref{nr3}), these limits exist.
As one easily checks :
$$\int_{B(0, \delta)} {{\varphi}_{\alpha}}^2 dv_{g_{\alpha}}=
{{\parallel u_{\alpha} \parallel  }_{\infty}}^{-2} {A_{\alpha}}^{-\frac{n}{2}} 
\int_{B_{x_{\alpha}}(\delta a_{\alpha})} {u_{\alpha}}^2 dv_g$$
and 
$$\int_{B(0, \delta)} {{\varphi}_{\alpha}}^{1+{\epsilon}_{\alpha}} dv_{g_{\alpha}}=
{{\parallel u_{\alpha} \parallel  }_{\infty}}^{-{(1+{\epsilon}_{\alpha})} }
{A_{\alpha}}^{-\frac{n}{2}} 
\int_{B_{x_{\alpha}}(\delta a_{\alpha})} {u_{\alpha}}^{1+{\epsilon}_{\alpha}} dv_g$$
$$=\left( {{\parallel u_{\alpha} \parallel  }_{\infty}}^{-{(1+{\epsilon}_{\alpha})} }
{A_{\alpha}}^{-\frac{n}{4}(1+{\epsilon}_{\alpha})} \right)
\left( {A_{\alpha}}^{-\frac{n}{4}(1+{\epsilon}_{\alpha})} 
\int_{B_{x_{\alpha}}(\delta a_{\alpha})} {u_{\alpha}}^{1+{\epsilon}_{\alpha}}
dv_g \right)
{A_{\alpha}}^{\frac{n}{2}{\epsilon}_{\alpha}}$$

\noindent Let first $\alpha$ goes to $0$ 
and then,  $\delta$ to $+\infty$.  By (\ref{az2})
and step 2, we have :

\begin{eqnarray} \label{R1}
\int_{\rR^n} {\varphi}^2 dv_{\xi}=
{C_0}^2
\end{eqnarray}
and 
\begin{eqnarray} \label{R2}
\int_{\rR^n} {\varphi} dv_{\xi}=
C_0 {\tilde{C}}_0^{\frac{n}{2}}
\end{eqnarray}

\noindent Now, let us compute  $\int_{\rR^n} {\mid \nabla \varphi
  \mid}_{\xi}^2 dv_{\xi}$. First, it is clear that  : 

\begin{eqnarray} \label{R3}
{\varphi}_{\alpha} \to \varphi \hbox{ in } C^1(B) \hbox{ as } \alpha
\to 0
\end{eqnarray}
for all compact ball $B$ in ${\rR}^n$.
Let ${\eta}_{\delta}(x)=\eta \left( {(2 \delta)}^{-1} \mid x \mid \right)$ where $\eta$ is
as in step 2.
Multiply $({\tilde E}_{\alpha})$ by ${\varphi}_{\alpha} {{\eta}_{\delta}}^2$ and
integrate over $\rR^n$. We check :
$$\int_{\rR^n} {< \nabla {\varphi}_{\alpha}, \nabla {\varphi}_{\alpha}
 {\eta}_{\delta}^2> }_{g_{\alpha}}
dv_{g_{\alpha}}
+\frac{2B_{\alpha}}{n {\parallel u_{\alpha}
\parallel}_{\infty}^{1-{\epsilon}_{\alpha}}}
\int_{\rR^n} {{\varphi}_{\alpha}}^{1+{\epsilon}_{\alpha}}
{{\eta}_{\delta}}^2
dv_{g_{\alpha}}
= \frac{k_{\alpha}}{2} \int_{\rR^n} {{\varphi}_{\alpha}}^2
{{\eta}_{\delta}}^2 dv_{g_{\alpha}}$$

\noindent Using (\ref{r2}), one  easily gets :
$$\lim_{\alpha \to 0} \frac{2B_{\alpha}}{n {\parallel u_{\alpha}
\parallel}_{\infty}^{1-{\epsilon}_{\alpha}}}=\frac{2}{n}{A_0(n)}^{-1}
C_0 \tilde{C}_0^{-\frac{n}{2}}$$
and then, by (\ref{r3}) and (\ref{R3}) :
$$\int_{\rR^n} {< \nabla \varphi ,\nabla \varphi 
{\eta}_{\delta}^2>}_{\xi} dv_{\xi}
+\frac{2}{n} {A_0(n)}^{-1} C_0 \tilde{C}_0^{-\frac{n}{2}}
\int_{\rR^n} {{\eta}_{\delta}}^2 \varphi dv_{\xi}$$
\begin{eqnarray} \label{R4}
= (1+\frac{2}{n}){A_0(n)}^{-1} \int_{\rR^n} {{\eta}_{\delta}}^2
{\varphi}^2 dv_{\xi}
\end{eqnarray}

\noindent We have 
$$\int_{\rR^n} {< \nabla \varphi ,\nabla \varphi 
{\eta}_{\delta}^2>}_{\xi} dv_{\xi}
=2\int_{\rR^n} {< \nabla \varphi, \nabla {\eta}_{\delta}>}_{\xi}
 \varphi {\eta}_{\delta} dv_{\xi}+\int_{\rR^n} {\mid \nabla \varphi 
\mid}_{\xi}^2 {{\eta}_{\delta}}^2 dv_{\xi}$$
$$\leq 2{\left(\int_{\rR^n} {\mid \nabla  {\eta}_{\delta}
\mid}_{\xi}^2 {\varphi}^2 dv_{\xi} \right)}^{\frac{1}{2}}
 {\left( \int_{\rR^n} {\mid \nabla \varphi 
\mid}_{\xi}^2 {{\eta}_{\delta}}^2 dv_{\xi} \right)}^{\frac{1}{2}}
+\int_{\rR^n} {\mid \nabla \varphi 
\mid}_{\xi}^2 {{\eta}_{\delta}}^2 dv_{\xi}$$

\noindent By (\ref{R1}) and since  ${\mid \nabla  
{\eta}_{\delta} \mid} \leq
\frac{\hbox{cst}}{{\delta}}$, one easily gets  :
\begin{eqnarray} \label{R5}
\lim_{\delta \to +\infty} \int_{\rR^n} {< \nabla \varphi, \nabla
  \varphi
  {\eta}_{\delta}^2>}_{\xi} dv_{\xi} = \int_{\rR^n} {\mid \nabla \varphi  
\mid}_{\xi}^2 dv_{\xi}
\end{eqnarray}

\noindent By (\ref{R1}), we know that $\varphi \in L^2(\rR^n)$. 
As a consequence, plugging (\ref{R5}) into (\ref{R4}) and using 
  (\ref{R2}), we have :

\begin{eqnarray} \label{R6}
 \int_{\rR^n} {\mid \nabla \varphi  
\mid}_{\xi}^2 dv_{\xi} =  {A_0(n)}^{-1} {C_0}^2
\end{eqnarray}

\noindent Now, let, for $u \in H_1^2(\rR^n)$:
$$I_{\xi}(u)=\frac{ \int_{\rR^n} {\mid \nabla u 
\mid}_{\xi}^2 dv_{\xi}
{(\int_{\rR^n} u dv_{\xi})}^{\frac{4}{n}}}
{{(\int_{\rR^n} u^2 dv_{\xi})}^{1+\frac{2}{n}}}$$

\noindent By the works of Carlen and Loss \cite{cl}, we know that :
$$\forall u \in H_1^2(\rR^n), I_{\xi}(u) \geq {A_0(n)}^{-1}$$

\noindent By (\ref{R1}), (\ref{R2}) and (\ref{R6}), we have :
$$I_{\xi}(\varphi)={A_0(n)}^{-1} \tilde{C}_0^2 $$
Since $\tilde{C}_0 \leq 1$,  it follows that ${\tilde{C}}_0=1$
( if ${\tilde{C}}_0<1$, we would have $I_{\xi}(\varphi)<{A_0(n)}^{-1}$ ). 
Therefore, $I_{\xi}(\varphi) ={A_0(n)}^{-1}$. 
Let $u$, $u \not\equiv 0$ and radially symetric, be an eigenfunction 
associated to $\lambda_1$, the first eigenvalue of the Laplacian on
the unit ball 
$\mathcal{B}$ in $\rR^n$ for radial functions with Neumann condition on
the boundary. Moreover, we may assume that u(0)=1. 
By Carlen and Loss \cite{cl}, we have :
$$\varphi=k v(\lambda x)$$
where $v(x)=u(x)-u(1)$. 
Now, by (\ref{R1}), (\ref{R2}) and since  ${\tilde{C}}_0=1$, we get :
$$\int_{\rR^n} {\varphi}^2 dv_{\xi}=
{\left(\int_{\rR^n} {\varphi} dv_g \right)}^2$$

\noindent We know that ( see theorem 1.3 in \cite{dhv} ) :
$$\int_{\rR^n} v^2 dv_{\xi}=\frac{n+2}{2} {u(1)}^2 \mid \mathcal{B} \mid$$
$$\int_{\rR^n} v dv_{\xi}=- \mid \mathcal{B} \mid u(1)$$

\noindent This gives then :
$${\lambda}^2={\lambda}_0^2$$
where 
$$\lambda_0^2={\left( \frac{n+2}{2} \right)}^{-\frac{2}{n}}
{ \mid \mathcal{B} \mid }^{\frac{2}{n}}$$

\noindent Let now :
\begin{eqnarray*}
r_{\nabla, \delta} & = &
\frac{\int_{B_{x_{\alpha}}(\delta a_{\alpha})} {\mid \nabla u_{\alpha} 
\mid}_g^2  R_{ij}(x_{\alpha})x^{i}x^{j}dv_g}
{\int_M {\mid \nabla u_{\alpha} {\eta}_{\alpha} \mid}^2 dv_g}\\
r_{1, \delta} & = & \frac{\int_{B_{x_{\alpha}}(\delta a_{\alpha})}
{( u_{\alpha} )}^{1+{\epsilon}_{\alpha}} 
R_{ij}(x_{\alpha})x^{i}x^{j}dv_g}
{\int_M {( u_{\alpha} {\eta}_{\alpha}) }^{1+{\epsilon}_{\alpha}} dv_g} \\
r_{2, \delta} & = & \frac{\int_{B_{x_{\alpha}}(\delta a_{\alpha})}
{( u_{\alpha} )}^2
R_{ij}(x_{\alpha})x^{i}x^{j}dv_g}
{\int_M {( u_{\alpha} {\eta}_{\alpha} )}^2 dv_g}  
\end{eqnarray*}
We recall that 
${\eta}_{\alpha}=\eta (c^{-1} r_{\alpha})$ where $c>0$ is small and
where
 $\eta$ is defined as before.
\noindent Using (\ref{az3}),
we easily see that
$$\lim_{\alpha \to 0}
\frac{\int_M {\mid \nabla u_{\alpha} {\eta}_{\alpha} \mid}_g^2 dv_g}
{\int_M
{\mid \nabla u_{\alpha}
 \mid}_g^2 dv_g}=1$$
We also get that, with step 2 and (\ref{az2}),  
$$\lim_{\alpha \to 0}
\frac{\int_M {( u_{\alpha} {\eta}_{\alpha} )}^2 dv_g}
{\int_M
u_{\alpha}^2 dv_g}=1$$
$$\lim_{\alpha \to 0}
\frac{\int_M {( u_{\alpha} {\eta}_{\alpha}) }^{1+{\epsilon}_{\alpha}} dv_g}
{\int_M
u_{\alpha}^{1+{\epsilon}_{\alpha}} dv_g}=1$$
\noindent Now, by an easy  proof by contradiction using step 2, 
(\ref{az3}) and (\ref{az2}), we see that   
$$\lim_{\delta \to \infty} \lim_{\alpha \to 0}
\frac{\int_M {\mid \nabla u_{\alpha}  \mid}_g^2 dv_g}
{\int_{B_{x_{\alpha}}(\delta a_{\alpha})}
{\mid \nabla u_{\alpha}
 \mid}_g^2 dv_g}=1$$

$$\lim_{\delta \to \infty} \lim_{\alpha \to 0}
\frac{\int_M u_{\alpha}^2 dv_g}
{\int_{B_{x_{\alpha}}(\delta a_{\alpha})}
u_{\alpha}^2 dv_g}=1$$
$$\lim_{\delta \to \infty} \lim_{\alpha \to 0}
\frac{\int_M u_{\alpha}^{1+{\epsilon}_{\alpha}} dv_g}
{\int_{B_{x_{\alpha}}(\delta a_{\alpha})}
u_{\alpha}^{1+{\epsilon}_{\alpha}} dv_g}=1$$
Here, $\lim_{\delta \to \infty} \lim_{\alpha \to 0}$
means that  $\alpha$ first goes to $0$ and then,
$\delta$ goes to  $+ \infty$. 
This implies that :
\begin{eqnarray*}
\lim_{\delta \to \infty} \lim_{\alpha \to 0} \frac{r_{\nabla,
 \delta}}{A_{\alpha}}
 & = &
\lim_{\delta \to \infty} \lim_{\alpha \to 0}
\frac{\int_{B_{x_{\alpha}}(\delta a_{\alpha})} {\mid \nabla u_{\alpha} 
\mid}_g^2  R_{ij}(x_{\alpha})x^{i}x^{j}dv_g}
{A_{\alpha} \int_{B_{x_{\alpha}}(\delta a_{\alpha})}
{\mid \nabla u_{\alpha}
 \mid}_g^2 dv_g}\\
\lim_{\delta \to \infty} \lim_{\alpha \to 0} \frac{r_{1,
 \delta}}{A_{\alpha}}
 & = & 
\lim_{\delta \to \infty} \lim_{\alpha \to 0}
\frac{\int_{B_{x_{\alpha}}(\delta a_{\alpha})}
{( u_{\alpha} )}^{1+{\epsilon}_{\alpha}} 
R_{ij}(x_{\alpha})x^{i}x^{j}dv_g}
{A_{\alpha} \int_{B_{x_{\alpha}}(\delta a_{\alpha})}
{( u_{\alpha}  )}^{1+{\epsilon}_{\alpha}}  dv_g} \\
\lim_{\delta \to \infty} \lim_{\alpha \to 0}
\frac{r_{2, \delta}}{A_{\alpha}} & = & 
\lim_{\delta \to \infty} \lim_{\alpha \to 0}
\frac{\int_{B_{x_{\alpha}}(\delta a_{\alpha})}
{( u_{\alpha} )}^2
R_{ij}(x_{\alpha})x^{i}x^{j}dv_g}
{A_{\alpha}\int_{B_{x_{\alpha}}(\delta a_{\alpha})}
{( u_{\alpha}  )}^2  dv_g}  
\end{eqnarray*}

\noindent Let $(y^1,..,y^n)$ be canonical coordinates in $\rR^n$ and
$(x^1,..,x^n)$ be exponential coordinates in $M$.
It is easy to see that, for a radial function $f$:
$$\int_{B(0, \delta)} f y^{i} y^{j} dv_{\xi}={\delta}^{ij} \frac{1}{n} 
\int_{B(0, \delta)} f {\mid y \mid}^2 dv_{\xi}$$
We also have :
$$\int_{B_{x_{\alpha}}(\delta a_{\alpha})} {u_{\alpha}}^p x^{i}x^{j} dv_{g}=
{\parallel u_{\alpha} \parallel}_{\infty}^{p} {A_{\alpha}}^{1+\frac{n}{2}}
\int_{B(0, {\delta})} {{\varphi}_{\alpha}}^p y^{i} y^{j} dv_{g_{\alpha}}$$
and :
$$\int_{B_{x_{\alpha}}(\delta a_{\alpha})}
{\mid \nabla u_{\alpha} \mid}_g^2 x^{i} x^{j} dv_g=
{\parallel u_{\alpha} \parallel}_{\infty}^2 {A_{\alpha}}^{\frac{n}{2}}
\int_{B(0, \delta)}
{\mid \nabla {\varphi}_{\alpha} \mid}_{g_{\alpha}} 
y^{i} y^{j} dv_{g_{\alpha}}$$

\noindent By these results and noting that
$\varphi$ is compactly supported, we have, for $\delta$ large enough :
$$\lim_{\alpha \to 0} \frac{r_{\nabla, \delta}}{A_{\alpha}}
=\frac{S_g(x_0)}{n}
\frac{\int_{\rR^n} {\mid \nabla \varphi \mid}_{\xi}^2 {\mid y \mid}^2
dv_{\xi}} {\int_{\rR^n} {\mid \nabla \varphi \mid}_{\xi}^2  dv_{\xi}}$$
$$\lim_{\alpha \to 0} \frac{r_{1, \delta}}{A_{\alpha}}
=\frac{S_g(x_0)}{n}
\frac{\int_{\rR^n}  \varphi 
{\mid y \mid}^2 dv_{\xi}}
{\int_{\rR^n} \varphi  dv_{\xi}}$$
$$\lim_{\alpha \to 0} \frac{r_{2, \delta}}{A_{\alpha}}
=\frac{S_g(x_0)}{n}
\frac{\int_{\rR^n} { \varphi }^2
{\mid y \mid}^2 dv_{\xi}}
{\int_{\rR^n}{ \varphi }^2 dv_{\xi}}$$
Then, for $\delta \geq {\lambda}_0$ :
$$\lim_{\alpha \to 0}
\frac{-\frac{1}{6} \left( -r_{\nabla, \delta}+(1+\frac{2}{n})r_{2, \delta}-
\frac{4}{n(1+{\epsilon}_{\alpha})} r_{1, \delta} \right)}
{A_{\alpha}}=\frac{{{\lambda}_0}^{-2} S_g(x_0)}{6n}
\Big( -\frac{\int_{\rR^n} {\mid \nabla v \mid}_{\xi}^2 {\mid y \mid}^2
dv_{\xi}} {\int_{\rR^n} {\mid \nabla v \mid}_{\xi}^2  dv_{\xi}}$$
$$+\frac{n+2}{n}\frac{\int_{\rR^n} { v }^2
{\mid y \mid}^2 dv_{\xi}}
{\int_{\rR^n}{ v }^2 dv_{\xi}}-\frac{4}{n(1+{\epsilon}_{\alpha})}
\frac{\int_{\rR^n}  v 
{\mid y \mid}^2 dv_{\xi}}
{ \int_{\rR^n} v   dv_{\xi}} \Big)$$

\noindent This expression has been computed in Druet, Hebey and Vaugon \cite{dhv}. 
We have :
$$\frac{-\frac{1}{6} \left( -r_{\nabla, \delta}+(1+\frac{2}{n})r_{2, \delta}-
\frac{4}{n(1+{\epsilon}_{\alpha})} r_{1, \delta} \right)}
{A_{\alpha}}=\frac{{\mid \mathcal{B}
\mid}^{-\frac{2}{n}}}{6n} \left( \frac{2}{n+2}+\frac{n-2}{{\lambda}_1} \right)
{\left(\frac{n+2}{2} \right)}^{\frac{2}{n}} S_g(x_{0}) $$

\noindent Hence, it is sufficient to prove that :

\begin{eqnarray} \label{R16}
\lim_{\delta \to \infty} \lim_{\alpha \to 0} 
\frac{r_{\nabla,\delta}-r_{\nabla}}{A_{\alpha}}=
\lim_{\delta \to \infty} \lim_{\alpha \to 0}
\frac{\int_{M-B_{x_{\alpha}}(\delta a_{\alpha})}
{\mid \nabla u_{\alpha}
{\eta}_{\alpha} \mid}_g^2  R_{ij}(x_{\alpha})x^{i}x^{j}dv_g}
{A_{\alpha} \int_M {\mid \nabla u_{\alpha} {\eta}_{\alpha} \mid}_g^2 dv_g}=0
\end{eqnarray}

\begin{eqnarray} \label{R17}
\lim_{\delta \to \infty} \lim_{\alpha \to 0} 
\frac{r_{1,\delta}-r_1}{A_{\alpha}}=
\lim_{\delta \to \infty} \lim_{\alpha \to 0}
\frac{\int_{M-B_{x_{\alpha}}(\delta a_{\alpha})}
{( u_{\alpha}
{\eta}_{\alpha} )}^{1+{\epsilon}_{\alpha}}  R_{ij}(x_{\alpha})x^{i}x^{j}dv_g}
{A_{\alpha}
\int_M {( u_{\alpha} {\eta}_{\alpha} )}^{1+{\epsilon}_{\alpha}}  dv_g}=0
\end{eqnarray} 

\begin{eqnarray} \label{R18}
\lim_{\delta \to \infty} \lim_{\alpha \to {\alpha}_0} 
\frac{r_{2,\delta}-r_{2}}{A_{\alpha}}=
\lim_{\delta \to \infty} \lim_{\alpha \to {\alpha}_0}
\frac{\int_{M-B_{x_{\alpha}}(\delta a_{\alpha})}
{( u_{\alpha}
{\eta}_{\alpha} )}^2 R_{ij}(x_{\alpha})x^{i}x^{j}dv_g}
{A_{\alpha} \int_M {( u_{\alpha} {\eta}_{\alpha} )}^2  dv_g}=0
\end{eqnarray} 
\noindent First, let us deal with (\ref{R18}).  Let :
$$T_{\alpha}= \left| \frac{\int_{M-B_{x_{\alpha}}(\delta a_{\alpha})}
{(\eta_{\alpha} u_{\alpha})}^2 R_{ij}(x_{\alpha})x^{i}x^{j}dv_g}
{A_{\alpha} \int_M {(\eta_{\alpha} u_{\alpha})}^2 dv_g} \right|$$
By (\ref{az2}) :

$$T_{\alpha}
\leq C
\frac{\int_{M-B_{x_{\alpha}}(\delta a_{\alpha})}
u_{\alpha}^2 {r_{\alpha}}^2 dv_g}
{A_{\alpha} }$$
Now, by step 3 :
$$T_{\alpha} \leq C
\frac{\int_{M-B_{x_{\alpha}}(\delta a_{\alpha})}
u_{\alpha}^{{\epsilon}_{\alpha}} {r_{\alpha}}^{2-n} 
{r_{\alpha}}^{\frac{n}{2}{\epsilon}_{\alpha}} dv_g}
{A_{\alpha}}$$
$$\leq C
\frac{\int_{M-B_{x_{\alpha}}(\delta a_{\alpha})}
u_{\alpha}^{{\epsilon}_{\alpha}} {r_{\alpha}}^{2-n} 
 dv_g}
{A_{\alpha} }
\leq C
\frac{{A_{\alpha}}^{1-\frac{n}{2}}
\int_{M-B_{x_{\alpha}}(\delta a_{\alpha})}
u_{\alpha}^{{\epsilon}_{\alpha}}  
 dv_g}
{A_{\alpha} }$$
To estimate this expression, we integrate  $(E_{\alpha})$ over
$M-B_{x_{\alpha}}(\delta a_{\alpha})$. We get :

\begin{eqnarray} \label{R19}
T_{\alpha} \leq C
\left( \frac{ {A_{\alpha}}^{-\frac{n}{2}}}{B_{\alpha}}
\int_{M-B_{x_{\alpha}}(\delta a_{\alpha})}
u_{\alpha} dv_g+\frac{{A_{\alpha}}^{1-\frac{n}{2}}}{B_{\alpha}}
\int_{\partial B_{x_{\alpha}}(\delta a_{\alpha})}
{\partial}_{\nu} u_{\alpha} d\sigma \right)
\end{eqnarray}

\noindent Let us prove that the second member of (\ref{R19}) goes to
0 if we let $\alpha$ goes to $0$ and
  ${\delta}$ to $\infty$. We have, using the definition of $A_\alpha$ :
  $$\frac{ {A_{\alpha}}^{-\frac{n}{2}}}{B_{\alpha}}    
\int_{M-B_{x_{\alpha}}(\delta a_{\alpha})}
u_{\alpha} dv_g
\leq \frac{ A_{\alpha}^{-\frac{n}{4}}}{B_{\alpha}}
{\left( \frac{\int_{M-B_{x_{\alpha}}(\delta a_{\alpha})}
u_{\alpha}^{1+{\epsilon}_{\alpha}} dv_g}
{\int_M u_{\alpha}^{1+{\epsilon}_{\alpha}} dv_g} 
\right)}^{\frac{1}{1+{\epsilon}_{\alpha}}}$$
By (\ref{r2}), we have :
$$\lim_{\alpha \to 0} \frac{ {A_{\alpha}}^{-\frac{n}{4}}}{B_{\alpha}}=C$$
Step 2 clearly implies that :
$$\lim_{\delta \to + \infty} \lim_{\alpha \to 0} {\left( \frac{\int_{M-B_{x_{\alpha}}(\delta a_{\alpha})}
u_{\alpha}^{1+{\epsilon}_{\alpha}} dv_g}
{\int_M u_{\alpha}^{1+{\epsilon}_{\alpha}} dv_g} 
\right)}^{\frac{1}{1+{\epsilon}_{\alpha}}}=0$$
Hence :
$$\lim_{\delta \to \infty} \lim_{\alpha \to 0}  
\frac{ {A_{\alpha}}^{-\frac{n}{2}}}{B_{\alpha}}
\int_{M-B_{x_{\alpha}}(\delta a_{\alpha})}
u_{\alpha} dv_g=0$$

\noindent  Now, if  $r_{\alpha}= \delta a_{\alpha}$, we
have :
$$\mid {\partial}_{\nu} u_{\alpha}(x) \mid
\leq \frac{{\parallel u_{\alpha} \parallel}_{\infty}}{{A_{\alpha}}^{\frac{1}{2}}}
{\parallel (\nabla \varphi)_g \parallel}_{L^{\infty}(\partial B(0, \delta))}$$

\noindent Since $\varphi$ is compactly supported ( see above ),
for $\delta$  large enough :
$${\parallel (\nabla {\varphi}_{\alpha})_{g_{\alpha}} \parallel}_{L^{\infty}(\partial
B(0, \delta))}
\to 0$$
Consequently, for $\delta$ large enough :
$$\lim_{\alpha \to 0} \frac{{A_{\alpha}}^{1-\frac{n}{2}}}{B_{\alpha}} 
\int_{\partial B_{x_{\alpha}}(\delta a_{\alpha})}
{\partial}_{\nu} u_{\alpha} d\sigma=0$$
By (\ref{R19}), this proves (\ref{R18}). 
To get (\ref{R16}) and (\ref{R17}), multiply 
$(E_{\alpha})$ by $\frac{{r_{\alpha}}^2
{{\eta}_{\alpha}}^2 u_{\alpha}}{A_{\alpha}}$ and integrate over
$M-B_{x_{\alpha}}(\delta a_{\alpha})$ :
$$-2  \int_{\partial B_{x_{\alpha}}(\delta a_{\alpha})}
(\partial_{\nu} u_{\alpha}) u_{\alpha} 
{r_{\alpha}}^2
{{\eta}_{\alpha}}^2 d\sigma
+2\int_{M-B_{x_{\alpha}}(\delta a_{\alpha})}
{\mid \nabla u_{\alpha} {\eta}_{\alpha} r_{\alpha} \mid}_g^2 dv_g
-2\int_{M-B_{x_{\alpha}}(\delta a_{\alpha})}
{\mid \nabla {\eta}_{\alpha} r_{\alpha} \mid}_g^2 {u_{\alpha}}^2 dv_g$$
\begin{eqnarray} \label{R20}
+\frac{4 B_{\alpha}}{n A_{\alpha}}
\int_{M-B_{x_{\alpha}}(\delta a_{\alpha})}
{u_{\alpha}}^{1+{\epsilon}_{\alpha}}
{r_{\alpha}}^2
{{\eta}_{\alpha}}^2 dv_g
=\frac{k_{\alpha}}{A_{\alpha}} 
\int_{M-B_{x_{\alpha}}(\delta a_{\alpha})}
{u_{\alpha}}^2
{r_{\alpha}}^2
{{\eta}_{\alpha}}^2 dv_g
\end{eqnarray}

\noindent As we did before, we use the fact that for $r_{\alpha}=\delta a_{\alpha}$ :
$$\mid {\partial}_{\nu} u_{\alpha}(x) \mid
\leq \frac{{\parallel u_{\alpha} \parallel}_{\infty}}{{A_{\alpha}}^{\frac{1}{2}}}
{\parallel (\nabla {\varphi}_{\alpha})_g \parallel}_{L^{\infty}(\partial
B(0, \delta))}$$ and :
$$u_{\alpha}(x)
\leq {\parallel u_{\alpha} \parallel}_{\infty}
{\parallel {\varphi}_{\alpha} \parallel}_{L^{\infty}(\partial B(0, \delta))}$$

\noindent This gives that for $\delta$ large enough, the boundary term goes to $0$. 
Moreover, it is clear that we have :
$$\int_{M-B_{x_{\alpha}}(\delta a_{\alpha})}
{\mid \nabla {{\eta}_{\alpha}} r_{\alpha} \mid}_g^2 u_{\alpha}^2
dv_g
\leq C \int_{M-B_{x_{\alpha}}(\delta a_{\alpha})}
 u_{\alpha}^2 
dv_g $$
By step 2, we obtain :
$$\lim_{\delta \to \infty} \lim_{\alpha \to 0} \int_{M-B_{x_{\alpha}}(\delta a_{\alpha})}
{\mid \nabla  r_{\alpha} {\eta}_{\alpha}^2 \mid}_g^2 {u_{\alpha}}^2  dv_g
=0$$
Observe that
the second member of (\ref{R20}) goes to $0$ when $\alpha \to 0$ and
$\delta \to \infty$. This easily follows from what we did when we proved
(\ref{R18}). Relation (\ref{R20}) then implies that :
\begin{eqnarray} \label{R21}
\lim_{\delta \to \infty} \lim_{\alpha \to 0} \int_{M-B_{x_{\alpha}}(\delta a_{\alpha})}
{\mid \nabla u_{\alpha} {\eta}_{\alpha} r_{\alpha} \mid}_g^2 dv_g = 0
\end{eqnarray}
and also that :
$$\lim_{\delta \to \infty}  \lim_{\alpha \to 0} 
\frac{4 B_{\alpha}}{n A_{\alpha}}
\int_{M-B_{x_{\alpha}}(\delta a_{\alpha})}
{u_{\alpha}}^{1+{\epsilon}_{\alpha}}
{r_{\alpha}}^2
{{\eta}_{\alpha}}^2 dv_g=0$$
which gives (\ref{R17}). In addition :
$$ \int_{M-B_{x_{\alpha}}(\delta a_{\alpha})}
{\mid \nabla u_{\alpha} {\eta}_{\alpha} r_{\alpha} \mid}_g^2 dv_g
= \int_{M-B_{x_{\alpha}}(\delta a_{\alpha})}
{\mid \nabla u_{\alpha} {\eta}_{\alpha} \mid}_g^2 {r_{\alpha}}^2 dv_g$$
$$+2 \int_{M-B_{x_{\alpha}}(\delta a_{\alpha})}
{<\nabla u_{\alpha} {\eta}_{\alpha}, \nabla r_{\alpha}>}_g 
u_{\alpha} {\eta}_{\alpha} r_{\alpha} dv_g
+\int_{M-B_{x_{\alpha}}(\delta a_{\alpha})}
{\mid \nabla r_{\alpha}  \mid}_g^2 {{\eta}_{\alpha} u_{\alpha}}^2 dv_g$$
For every $x, y, \epsilon > 0$, we have : $xy \leq \frac{1}{2}( \epsilon x^2
+\frac{1}{\epsilon}  y^2)$. Noting that :
$$ \int_{M-B_{x_{\alpha}}(\delta a_{\alpha})}
{<\nabla u_{\alpha} {\eta}_{\alpha}, \nabla r_{\alpha}>}_g 
u_{\alpha} {\eta}_{\alpha} r_{\alpha} dv_g$$
$$\geq - {\left( \int_{M-B_{x_{\alpha}}(\delta a_{\alpha})}
{\mid \nabla u_{\alpha} {\eta}_{\alpha} \mid}_g^2 {r_{\alpha}}^2 dv_g
\right)}^{\frac{1}{2}}
{\left( \int_{M-B_{x_{\alpha}}(\delta a_{\alpha})}
{\mid \nabla r_{\alpha}  \mid}_g^2 {{\eta}_{\alpha} u_{\alpha}}^2 dv_g
\right)}^{\frac{1}{2}}$$
we get :
$$\int_{M-B_{x_{\alpha}}(\delta a_{\alpha})}
{\mid \nabla u_{\alpha} {\eta}_{\alpha} r_{\alpha} \mid}_g^2 dv_g
\geq (1-\epsilon)
\int_{M-B_{x_{\alpha}}(\delta a_{\alpha})}
{\mid \nabla u_{\alpha} {\eta}_{\alpha} \mid}_g^2 {r_{\alpha}}^2 dv_g$$
$$+(1-\frac{1}{\epsilon})
\int_{M-B_{x_{\alpha}}(\delta a_{\alpha})}
{\mid \nabla r_{\alpha}  \mid}_g^2 {({\eta}_{\alpha} u_{\alpha})}^2 dv_g$$

\noindent Using (\ref{R21}) and the fact that 
$\lim A_{\alpha} \int_M {\mid \nabla u_{\alpha} \eta_{\alpha}
    \mid}_g^2 dv_g= A_0(n)^{-1}$, we then clearly get (\ref{R16}). 
Finally,  this proves step 5.

\begin{step}
We prove the theorem.
\end{step}

\noindent Let, for $u \in H_1^2(M)$ :
$$I_{g,\alpha} (u)=I_{\alpha}(u)-(\alpha_0-\alpha) 
{(\int_M {\mid u \mid}^{1+{\epsilon}_{\alpha}}
dv_g)}^{\frac{4}{n(1+{\epsilon}_{\alpha})}}$$

\noindent $\bf{a-}$
$We$ $first$ $prove$ $that$ :
\begin{eqnarray} \label{R7}
\lim_{\alpha \to 0} \frac{{A_0(n)}^{-1} - I_{g,\alpha}({\eta}_{\alpha} u_{\alpha})}
{A_{\alpha}}={\alpha}_0
\end{eqnarray}
By (\ref{r31}), (\ref{r32}) and (\ref{r33}), one can check that :
\begin{eqnarray} \label{R8}
\lim_{\alpha \to 0} 
\frac{I_{g,\alpha}(u_{\alpha})-
I_{g,\alpha}({\eta}_{\alpha} u_{\alpha})}{A_{\alpha}}=0
\end{eqnarray}
Moreover, we have :
$$I_{g,\alpha}(u_{\alpha})=I_{\alpha}(u_{\alpha})-(\alpha_0-\alpha)A_{\alpha}$$
Since $\alpha \to 0$ and $I_{\alpha}(u_{\alpha}) 
\leq {A_0(n)}^{-1}$, we get :
\begin{eqnarray} \label{R9}
\liminf_{\alpha \to 0} \frac{{A_0(n)}^{-1} - I_{g,\alpha}({\eta}_{\alpha} u_{\alpha})}
{A_{\alpha}} \geq {\alpha}_0
\end{eqnarray}
In addition, we can also write, by (\ref{R8})
\begin{eqnarray*} 
\limsup_{\alpha \to 0} 
\frac{{A_0(n)}^{-1} - I_{g,\alpha}({\eta}_{\alpha} u_{\alpha})}
{A_{\alpha}}=\limsup_{\alpha \to 0}
\frac{{A_0(n)}^{-1} - I_0( u_{\alpha})+\alpha_0 A_{\alpha}}
{A_{\alpha}}
\end{eqnarray*}
By definition of $\alpha_0$, we have $ I_0( u_{\alpha}) \geq \mu_0 = 
A_0(n)^{-1}$. This implies that :
\begin{eqnarray} \label{R10}
\limsup_{\alpha \to 0} \frac{{A_0(n)}^{-1} - I_{g,\alpha}({\eta}_{\alpha} u_{\alpha})}
{A_{\alpha}} \leq \alpha_0
\end{eqnarray}
(\ref{R7}) then comes from (\ref{R8}), (\ref{R9}) and (\ref{R10}). \\

\noindent $\bf{b-}$
$We$ $prove$ $that$ :
\begin{eqnarray} \label{R23}
\int_M {\mid \nabla {\eta}_{\alpha} u_{\alpha}
\mid}_{\xi}^2 dv_{\xi} -
\int_M {\mid \nabla {\eta}_{\alpha} u_{\alpha} \mid}_g^2 dv_g 
= -\frac{1}{6}
\int_M {\mid \nabla {\eta}_{\alpha} u_{\alpha} \mid}_{\xi}^2 
R_{ij}(x_{\alpha}) x^{i}x^{j} dv_g
+O(1)
\end{eqnarray}

\noindent First note that the limit of right-hand side member of 
(\ref{R23}) exists.
We have 

\begin{eqnarray} \label{R24}
\int_M {\mid \nabla {\eta}_{\alpha} u_{\alpha} \mid}_g^2 dv_g=
\int_M {\mid \nabla {\eta}_{\alpha} u_{\alpha} \mid}_{\xi}^2 dv_g+
\int_M (g^{ij}-{\delta}^{ij}){\partial}_i u_{\alpha} {\partial}_j u_{\alpha}
{\eta}_{\alpha}^2 dv_g+C_1(\alpha)
\end{eqnarray}
where $C_1(\alpha)$ stands for the terms in which the derivatives of 
$\eta_{\alpha}$ appear. 
Since  $supp(\nabla {\eta}_{\alpha}) \subset
M-B_{x_{\alpha}}(\frac{c}{2})$ and by step  2, (\ref{az3}) and
(\ref{az2}),
we see that  $C_1(\alpha) \to 0$
when $\alpha \to 0$. 
We write that, for $\delta >0$, 
$$\left| 
\int_M (g^{ij}-{\delta}^{ij}){\partial}_i u_{\alpha} {\partial}_j u_{\alpha}
{\eta}_{\alpha}^2 dv_g \right| \leq \left| 
\int_{{B_{x_{\alpha}}(\delta a_{\alpha})}}  (
g^{ij}-{\delta}^{ij}){\partial}_i u_{\alpha} {\partial}_j u_{\alpha}
 dv_g \right|$$
$$+\left| 
\int_{M-{B_{x_{\alpha}}(\delta a_{\alpha})}}  (g^{ij}-{\delta}^{ij}){\partial}_i u_{\alpha} {\partial}_j u_{\alpha}
{\eta}_{\alpha}^2 dv_g \right|$$
Using the Cartan Hadamard expansion of the metric  $g$, we get that
$$\left| 
\int_M (g^{ij}-{\delta}^{ij}){\partial}_i u_{\alpha} {\partial}_j u_{\alpha}
{\eta}_{\alpha}^2 dv_g \right|
\leq C  \left|
\int_{B_{x_{\alpha}}(\delta a_{\alpha})} 
R^i{}_{kl}{}^j(x_{\alpha}) \partial_i u_{\alpha} \partial_j u_{\alpha}
x^k x^l 
dv_g \right| $$
$$+C \int_{B_{x_{\alpha}}(\delta a_{\alpha})}  
{\mid \nabla u_{\alpha}
\mid}_g^2 r_{\alpha}^3 dv_g
+C \int_{M-B_{x_{\alpha}}(\delta a_{\alpha})} 
{\mid \nabla u_{\alpha}
\mid}_g^2 r_{\alpha}^2 dv_g$$
where $(R^i{}_{kl}{}^j(x_{\alpha}))$ are the components of the Riemann
curvature of $g$ in exponential map at $x_{\alpha}$.
One gets from (\ref{R21}) that the third term of this expression is small if $\delta$ is
large. The second term goes to $0$ when 
$\alpha$ tends to $0$. It can be seen by writing that, on 
$B_{x_{\alpha}}(\delta a_{\alpha})$, $r_{\alpha} \leq \delta
  a_{\alpha}$. 
We now prove that the first term goes to $0$ with   $\alpha$. We write
that 
$$\left|
\int_{B_{x_{\alpha}}(\delta a_{\alpha})} 
R^i{}_{kl}{}^j(x_{\alpha}) \partial_i u_{\alpha} \partial_j u_{\alpha}
x^k x^l  
dv_g \right|$$
$$\leq C 
{\parallel u_{\alpha} \parallel}_{\infty}^2 A_{\alpha}^{\frac{n}{2}}
\left| \int_{B(0,\delta)} R^i{}_{kl}{}^j(x_{\alpha}) 
\partial_i \varphi_{\alpha} \partial_j \varphi_{\alpha}
x^k x^l  dv_{g_{\alpha}} \right|$$
where $\varphi$ is defined as in step 1. Now, since ${\varphi}_{\alpha} 
\to {\varphi}$ in  
$C^1(B(0,\delta))$ when $\alpha \to 0$ and  since  $\varphi$ is
radially 
symmetric, we get that
$$\lim_{\alpha \to 0} R^i{}_{kl}{}^j(x_{\alpha}) \partial_i u_{\alpha} \partial_j u_{\alpha}
x^k x^l  =0$$
Together with  (\ref{r11}), this proves that, for all  $\delta$,  
$$\lim_{\alpha \to 0} 
\int_{B_{x_{\alpha}}(\delta a_{\alpha})} 
R^i{}_{kl}{}^j(x_{\alpha}) \partial_i \varphi_{\alpha} \partial_j 
\varphi_{\alpha}
x^k x^l  
dv_g  =0$$
We finally obtain that 
\begin{eqnarray} \label{R25}
\lim_{\alpha \to 0} \int_M (g^{ij}-{\delta}^{ij}){\partial}_i u_{\alpha} {\partial}_j u_{\alpha}
{\eta}_{\alpha}^2 dv_g=0
\end{eqnarray}
To conclude, we write that, by the 
Cartan Hadamard expansion of $g$,  
\begin{eqnarray} \label{R26}
\int_M {\mid \nabla {\eta}_{\alpha} u_{\alpha} \mid}_{\xi}^2 dv_g=
\int_M {\mid \nabla {\eta}_{\alpha} u_{\alpha} \mid}_{\xi}^2 dv_{\xi}
+\frac{1}{6}
\int_M {\mid \nabla {\eta}_{\alpha} u_{\alpha} \mid}_{\xi}^2 
R_{ij}(x_{\alpha}) x^{i}x^{j} dv_g+O(1)
\end{eqnarray}

\noindent We then get  (\ref{R23}) from (\ref{R24}), (\ref{R25}) and 
(\ref{R26}). \\

\noindent $\bf{c-}$
$We$ $prove$ $that$ : 
\begin{eqnarray} \label{R22}
\lim_{\alpha \to 0}
\frac{I_{\xi,\alpha}({\eta}_{\alpha}u_{\alpha})-I_{g,\alpha}({\eta}_{\alpha}
u_{\alpha})}{A_{\alpha}} 
={A_0(n)}^{-1} \frac{{\mid \mathcal{B} \mid}^{-\frac{2}{n}}}{6n}
\left( \frac{2}{n+2}+\frac{n-2}{{\lambda}_1} \right)
{\left(\frac{n+2}{2}\right)}^{\frac{2}{n}} S_g(x_0)
\end{eqnarray}
\noindent $where$ $I_{\xi}$ $is$ $defined$ $as$ $above$. 

\noindent Let :
\begin{eqnarray*}
t_1 & = & \frac{\int_M (\eta_{\alpha} u_{\alpha})^{1+\epsilon_{\alpha}}
dv_{\xi}-\int_M (\eta_{\alpha} u_{\alpha})^{1+\epsilon_{\alpha}}
dv_g}{ \int_M (\eta_{\alpha} u_{\alpha})^{1+\epsilon_{\alpha}}
dv_g}\\
t_2 & = & \frac{\int_M (\eta_{\alpha} u_{\alpha})^2
dv_{\xi}-\int_M (\eta_{\alpha} u_{\alpha})^2 dv_g}
{ \int_M (\eta_{\alpha} u_{\alpha})^2
dv_g}\\
t_{\nabla} & = &  \frac{\int_M {\mid \nabla \eta_{\alpha} 
u_{\alpha} \mid}_{\xi}^2
dv_{\xi}-\int_M {\mid \nabla \eta_{\alpha} 
u_{\alpha} \mid}_g^2
dv_g}
{\int_M {\mid \nabla \eta_{\alpha} 
u_{\alpha} \mid}_g^2
dv_g }
\end{eqnarray*}

\noindent By the Cartan Hadamard expansion of $g$, we have :
$$dv_{\xi}=\left(1+\frac{1}{6} R_{i,j}(x_{\alpha}) x^i
  x^j+O(r_{\alpha}^3)
\right) dv_g$$
Coming back to the notations of step 5, we then get :
\begin{eqnarray} \label{ab1}
\lim_{\alpha \to 0} \frac{t_1}{A_{\alpha}}=\lim_{\alpha \to 0} \frac{1}{6}\frac{r_1}{A_{\alpha}}
\end{eqnarray}
and :
\begin{eqnarray} \label{ab2}
\lim_{\alpha \to 0} \frac{t_2}{A_{\alpha}}=\lim_{\alpha \to 0} \frac{1}{6}\frac{r_2}{A_{\alpha}}
\end{eqnarray}
From (\ref{R23}), we also have :
\begin{eqnarray} \label{ab3}
\lim_{\alpha \to 0} \frac{t_{\nabla}}{A_{\alpha}}=\lim_{\alpha \to 0} \frac{1}{6}\frac{r_{\nabla}}{A_{\alpha}}
\end{eqnarray}

\noindent We write :
$$I_{\xi,\alpha}(u_{\alpha} {\eta}_{\alpha})-I_{g,\alpha}
(u_{\alpha} {\eta}_{\alpha})=
I_{g,\alpha}(u_{\alpha}{\eta}_{\alpha})
\frac{(1+t_{\nabla} ){(1+t_1)}
^{\frac{4}{n(1+{\epsilon}_{\alpha})}}}
{{(1+t_2)}^{1+\frac{2}{n}}}-I_{g,\alpha}(u_{\alpha}
{\eta}_{\alpha})$$
(\ref{R22}) then follows by (\ref{s1}),
(\ref{ab1}),(\ref{ab2}), (\ref{ab3}) and the fact that $\lim_{\alpha
\to 0} I_{g,{\alpha}}(u_{\alpha} \eta_{\alpha})={A_0(n)}^{-1}$.\\

\noindent $\bf{d-}$
$Conclusion$

\noindent By H$\ddot{o}$lder's inequality and Carlen and
Loss \cite{cl}, we have :
$$I_{\xi,\alpha}({\eta}_{\alpha}u_{\alpha}) \geq 
\frac{\int_M {\mid \nabla {\eta}_{\alpha} u_{\alpha} \mid}_{\xi}^2  dv_{\xi} 
{\left( \int_M {\eta}_{\alpha} u_{\alpha} dv_{\xi}
  \right)}^{\frac{4}{n}}}{\int_M {({\eta}_{\alpha}u_{\alpha})}^2
dv_{\xi}} 
\geq {A_0(n)}^{-1}$$
We have then :
$$I_{\xi,\alpha}({\eta}_{\alpha}u_{\alpha}) -
I_{g,\alpha}({\eta}_{\alpha}u_{\alpha}) \geq
{A_0(n)}^{-1}-I_{g,\alpha}g({\eta}_{\alpha}u_{\alpha})$$
Dividing this inequality by $A_{\alpha}$ and recalling that $B_0={\alpha}_0 A_0(n)$,
we get from (\ref{R7}) and (\ref{R22}) that :
$$B_0 \leq \frac{{\mid \mathcal{B} \mid}^{-\frac{2}{n}}}{6n}
\left( \frac{2}{n+2}+\frac{n-2}{{\lambda}_1} \right)
{\left(\frac{n+2}{2}\right)}^{\frac{2}{n}} S_g(x_0)$$
and then :
$$B_0 \leq \frac{{\mid \mathcal{B} \mid}^{-\frac{2}{n}}}{6n}
\left( \frac{2}{n+2}+\frac{n-2}{{\lambda}_1} \right)
{\left(\frac{n+2}{2}\right)}^{\frac{2}{n}} \max_{x \in M}{S_g(x)}$$
This ends the proof of the theorem.

%%%%%%%%%%%%%%%%%%%%%%%%%%%%%%%%%%%%%%%%%%%%%%%%%%%%%%%%%%%%%%%%%%%%%%%%%%%%%


\begin{thebibliography}{9}


\bibitem{a1}
{\sc T.Aubin}--
Some Nonlinear problems in Riemannian Geometry,     
{\em Berlin Springer-Verlag, 1998. }


\bibitem{al}
{\sc T.Aubin} {\sc Y.Y.Li}--
On the best Sobolev inequality,
{\em Journal de Math\'ematiques Pures et Appliqu\'ees,
78, No4, 1999, p.353-382. } 

\bibitem{cl}
{\sc E.A.Carlen} {\sc M.Loss}--
Sharp constant in Nash's inequality,  {\em International Mathematics
Notices, 7, 1993, p.213-215. }
 

\bibitem{dd}
{\sc Z.Djadli} {\sc O.Druet }--  
Extremal functions for optimal Sobolev inequalities on compact
manifolds,  
{\em Calculus of Variations and Partial Differential Equations, 12, No1, 2001,
p.59-84.}
 

\bibitem{d}
{\sc O.Druet}--  
The best constant problem in Sololev inequalities, 
{\em Mathematische
Annalen, 314, No2, 1999, p.327-346. }



\bibitem{dhv}
{\sc O.Druet} {\sc E.Hebey} {\sc M.Vaugon}--
Optimal Nash's inequalities on riemannian manifolds,   {\em 
International Mathematics Research Notices, 14, 1999, p.735-779. }



\bibitem{hv}
{\sc E.Hebey} {\sc M.Vaugon}--
Meilleures constantes dans le th\'eor\`eme d'inclusion de Sobolev,  
{\em Annales de l'Institut Henri Poincar\'e, Analyse non-lin\'eaire, vol. 13, 1996,
57-93. }

\bibitem{hv1}
{\sc E.Hebey} {\sc M.Vaugon}--
From best constants to critical functions,  
{\em Mathematische Zeitschrift, 237, No 4, 2001,
737-767. }


\bibitem {heb}
{\sc E.Hebey}--
Nonlinear Analysis on Manifolds : Sobolev Spaces and Inequalities,
{\em Lecture Notes, Courant Institute, Vol. 5, 1999. }
 

\bibitem{h}
{\sc E.Humbert}--
Best constants in the $L^2-$Nash inequality,  {\em to appear in
  Proceedings of the Royal Society of Edinburgh.}


\end{thebibliography}
\end{document}